\documentclass[fleqn,11pt]{article}
\usepackage{amsmath}
\usepackage{color}
\usepackage{amssymb}
\usepackage{graphicx}

\renewcommand{\baselinestretch}{1.3}
\newtheorem{theorem}{\bf Theorem}[section]

\newtheorem{remark}{\bf Remark}

\newtheorem{exam}[theorem]{\bf Example }

    {\par \noindent {\it Proof. }}%
    {\par \indent}
    {\par \noindent {\bf Counterexample:}}%
    {\par \indent}

\def\E{{\rm E}}
\def\var{{\rm Var}}

\newcommand{\be}{\begin{equation}}
\newcommand{\bes}{\begin{displaymath}}
\newcommand{\ee}{\end{equation}}
\newcommand{\ees}{\end{displaymath}}

\usepackage{latexsym}% necessary for end of proof symbol according to
\usepackage{rotate,epsfig}% necessary for figures

\begin{document}

\title{ Bayesian prediction of minimal repair times of a series system based on hybrid censored sample of components' lifetimes  under Rayleigh distribution }
\author{ {{ S. M. T. K. MirMostafaee}}${}^{1, }$\footnote{Corresponding author.\newline E-mail addresses:
{\it m.mirmostafaee@umz.ac.ir} (S. M. T. K. MirMostafaee), {\it morteza.amini@ut.ac.ir}
(Morteza Amini), {\it a.asgharzadeh@umz.ac.ir} (A. Asgharzadeh).}\ , {{ Morteza Amini}}${}^{2}$ and {{{ A. Asgharzadeh}}${}^{1}$ \vspace{2mm}}\\
{\small {\it $~^{1}$Department of Statistics, Faculty of Mathematical Sciences,
  University of Mazandaran, }}\vspace{-0.1cm}\\
{\small {\it P. O. Box 47416-1467,  Babolsar,\ Iran}}\\
{\small {\it $~^{2}$Department of Statistics, School of
Mathematics, Statistics and Computer }}\vspace{-0.1cm}\\ {\small {\it  Sciences, College of Science, University of Tehran, P.O. Box 14155-6455,  Tehran, Iran}}}
%\date{}
\maketitle
\begin{abstract}
In this paper, we develop Bayesian predictive inferential procedures for
prediction of repair times of a series system, applying
a minimal repair strategy, using the information contained in
an independent observed hybrid censored sample of the lifetimes
of the components of the system, assuming the underlying distribution of the lifetimes to be Rayleigh distribution. An illustrative
real data example and a simulation study are presented for the purpose of illustration and comparison of the proposed predictors.
\end{abstract}


\section{ Introduction}

Suppose that a series system with $k$  repairable and identical components is
under operation. We assume that these $k$ components work
independently and that the lifetimes of all components are
identically distributed. The lifetime of this
system is equal to the minimum lifetime of its components. We also assume that the system is repaired
using the minimal repair strategy. In a minimal repair strategy, after each
failure, the corrupt component is immediately repaired and restored to
its original condition, thus putting back the
system into operation. Therefore, the state of the
system after a repair is the same as it was immediately before
corresponding failure. We shall assume that the time needed for repair and
replacement is negligible. The minimal repair
times possess the same joint distribution as upper record values
from the distribution of the lifetime of the system, that is the
distribution of record values from the distribution of the minimum of a sample of size $k$ (see Barlow
and Hunter, 1960). It has been verified that the
sequence of record values, from the distribution of minimum in a
sample of size $k$, and the sequence of
$k$-record values, from the parent distribution, are identically distributed (see Arnold et. al., 1998, p. 43).

The results of this paper focus on predicting the minimal repair times of a series system based on an available hybrid censored sample of its components' lifetimes.  Consider a sample
of $n$ units placed on a life-test at time 0. In  Type-I censoring
scheme, the test is terminated at a pre-fixed time $T$, while in
Type-II censoring scheme, it is terminated as soon as
a pre-determined number, $r$ ($r\leq n)$, of units
fail. Under Type-I censoring scheme, the duration of the life-test is guaranteed, while the number of failures is random, which might result in a low efficiency, when the number of failures is small. In Type-II censoring scheme, the level of efficiency is guaranteed, as the number of failures, $r$, is pre-fixed, while the exact time of the $r^{\textrm{th}}$ failure is random, thus the duration of the experiment may
end up being too long. The mixture of the Type-I and Type-II censoring schemes, called {\em hybrid censoring} scheme, proposed first by Epstein (1954), reduces the mentioned disadvantages.
Under a Type-I hybrid censoring scheme, the experiment is terminated as soon as either $r$ units fail or
the time $T$ is reached.

Hybrid censoring has received a
considerable attention in the context of
reliability and life-testing. Many authors, including Draper and Guttman (1987), Kundu and Gupta (1988),
Ebrahimi (1992), Childs et al. (2003), Kundu (2007) and Kundu and
Banerjee (2008), have developed statistical inference based on hybrid censored sample.
For a comprehensive review of hybrid censoring, see Balakrishnan and Kundu (2012).

The real data used in this paper includes the number of revolutions to failure of ball bearings under a life test, accelerated by hybrid censoring. The ball bearings are identical, thus the components' lifetimes follow the same distribution. The test is performed before placement of $k$ identical ball bearings in a machine. The machine, made up of the $k$ such identical components fails as soon as the first ball bearing fails, that is that the machine is a series system of $k$ identical components. In the case of the failure, the physical and statistical (black box) minimal repair of the system are equivalent and are performed by minimal repair of the failed component. Our aim here is to predict the minimal repair times of this machine, using the information achieved from the available censored sample.

In this paper, we assume that the underlying lifetime
distribution is the two parameter Rayleigh distribution, with cumulative distribution function (cdf),
\begin{equation}\label{eq:cdfre}
F(x;\mu,\sigma)=\left\{
\begin{array}{l l}
0, & x\leq\mu,\\
1- \exp\left\{-\frac{(x-\mu)^2}{2\sigma}\right\}, & x>\mu,
\end{array}\right.
\end{equation}
where $\mu\in\mathbb{R}$ and $\sigma>0.$ When $\mu=0$, the distribution \eqref{eq:cdfre} is called the \textit{scaled Rayleigh distribution}. The corresponding probability density function (pdf) of \eqref{eq:cdfre} is
\begin{eqnarray}\label{eq:pdfre}
f(x;\mu,\sigma)=\frac{(x-\mu)}{\sigma} \exp\left\{-\frac{(x-\mu)^2}{2\sigma}\right\}, \ \  x>\mu.
 \end{eqnarray}

The Rayleigh distribution is widely applied in several areas of statistics, partly because of its linear and increasing failure rate, which makes it an appropriate distribution for modeling the lifetime distribution of components, which age rapidly with time. Several types of electro-vacuum devices have this feature (Polovko, 1968). The Rayleigh distribution was originally introduced by Lord Rayleigh (1880, 1919) in connection with a problem in the field of acoustics. Wide applications of the Rayleigh distribution in lifetime analysis is mentioned by many authors, including Polovko (1968), Johnson et al. (1994), Dyer and Whisenand (1973), Gross and Clark (1976),  Balakrishnan (1989) and Lawless (2003). The Rayleigh distribution relates to a number of well-known life distributions such as generalized extreme value, Weibull and Chi-square distributions (see Dey and Dey, 2014). There are many papers dealing with estimation and prediction under Rayleigh distribution, including Howlader (1985), Howlader and Hossain (1995), Fern\'andez (2000), Raqab and Madi (2002), Ali Mousa and Al-Sagheer (2005), Dey and Das (2007), Soliman and Al-Aboud (2008), Khan et al. (2010), Dey and Dey (2012) and Abdel-Hamid and Al-Hussaini (2014).

Considerable work has been done on prediction of future records and order statistics, using parametric and
nonparametric inferential methods. The following papers consider the Rayleigh distribution as the underlying distribution and develop predictive inferential methods for records and order statistics.
Howlader (1985) obtained the highest posterior density (HPD) prediction intervals for future order statistics from an independent sample, based on an observed sample of order statistics. Fernandez (2000) considered the problem of Bayesian prediction of a future observation based on an observed Type-II doubly censored sample. Raqab and Madi (2002) developed Bayesian prediction of the total time on a test based on doubly censored sample.
Ali Mousa and Al-Sagheer (2005) considered Bayesian prediction of a progressive censored sample on the basis of an observed  progressively Type-II censored sample. Soliman and Al-Aboud (2008) derived Bayesian prediction intervals for a future record value based on an observed sample of record values. Recently, Khan et al. (2010) have develpoed Bayesian inference about a future order statistic on the basis of a doubly censored sample.

The observed sample and the predicted future
observation might be either of the same type or of different types.
Recently, Ahmadi and MirMostafaee (2009), Ahmadi and Balakrishnan (2010), Ahmadi et al. (2010) and MirMostafaee and Ahmadi (2011), have considered the problem of predicting future
records from a sequence of observations on the basis of order
statistics observed from another independent sample and vice
versa. According to our best knowledge, there is a few number of similar works in the context of
prediction of record values based on an available censored sample. These works use simple Type-II censored sample as the available data. The hybrid censoring is a mixture of the Type-I and Type-II censoring schemes.
Thus the termination time of the experiment in hybrid censoring would decay stochastically with respect to
Type-I and Type-II censoring schemes.

In this paper, we obtain several Bayesian point predictors as
well as Bayesian prediction intervals for a future
repair time of a $k$-component series system, applying a minimal
repair strategy, on the basis of observed hybrid censored sample
of its components lifetime, when the underlying distribution is
assumed to be scaled or two parameter Rayleigh with cdf
\eqref{eq:cdfre}, and $\mu=0$ or $\mu\in\mathbb{R}$, respectively.

The rest of this paper is organized as follows. Section 2 is devoted to the description of the model.
In Section 3, we develop the main results. Sections 4 contains a real data example, which
illustrate the results. A simulation study is presented in Section 5 for comparison of the performance of the proposed predictors.

\section{The model}

Suppose that $n$ identical components are under a life test. Let $X_1,\ldots,X_n$ denote the lifetimes of the
experimental units and $X_{1:n} < \ldots< X_{n:n}$ stand for the
corresponding order statistics. Furthermore, suppose that the experiment is terminated according to a hybrid censoring strategy. Under a hybrid
censoring scheme, the experiment is terminated at
time $T_0=\min\{X_{r:n},T\}$, $r\leq n$, where $r$ and $T$ are pre-determined values. The
observed hybrid censored sample is therefore ${\bf X}=(X_{1:n}
\ldots, X_{d:n})$, where { $d=\max\{s: s\leq r,\; X_{s:n}\leq T\}$}.
To simplify the notation, we henceforth denote the hybrid
censored sample by $(X_1,X_2,\cdots,X_d)$.

When $X_1,\ldots,X_n$ are independent and
identically distributed (iid) with the common absolutely continuous cdf
$F$ and corresponding pdf $f$, the joint pdf of the hybrid
censored sample is
\begin{equation}\label{hyb}
f_{X_1,\ldots,X_d}(x_1,\ldots,x_d)=C\prod_{i=1}^{d}f(x_i)[1-F(T_0)]^{n-d},
\end{equation}
where $C$ is the normalizing constant.

Suppose that the test is performed before placement of $k$ identical components in a machine. Assume further that the machine, made up of the $k$ identical components, fails as soon as the first component fails, that is the machine is a series system of $k$ identical components.

In the sequel, we develop Bayesian prediction of the repair times of this machine, repaired using a minimal repairing strategy, based on the available hybrid censored sample $(X_1,X_2,\cdots,X_d)$. The future $m^{\rm th}$ repair
time of a series system, with $k$ independent and identical components, the lifetimes of
which have the same distribution as $X_1$, is denoted by $U_{m(k)}$. 

Let the sequence $\{Y_i\}_{i=1}^{\infty}$ be independent and identically distributed random
variables, independent of the sample of component lifetimes,
$X_1,\ldots,X_n$, each of which having the same distribution as $X_1$. Also, let $Y_{j:n}$ 
stand for the $j^{\rm th}$ order statistic among $Y_1,\ldots,Y_n$. The repair time, $U_{m(k)}$, is distributed as
$(T({m,k})-k+1)^{\rm th}$ order statistic among $Y_1,\ldots,Y_{T({m,k})}$, where $T({1,k})=k$
and 
$$T({m, k})=\min\{j:j>T({m-1, k}), ~Y_j>Y_{T({m-1, k})-k+1:T({m-1, k})}\},$$
for  $m\geq 2$.

Indeed, $U_{1(k)},U_{2(k)},\ldots$ is the sequence of current
$k^{\mbox{th}}$ largest $Y$s yet seen (see Arnold et al., 1998),
that is $U_{m(k)}$ is identically distributed to the $m^{\rm th}$
$k$-record value.

The sequence $\{U_{m(k)}\}_{m\geq1}$ from a cdf $F$ is identical
in distribution to a record sequence $\{U_{m(1)}\}_{m\geq1}$ from
the cdf of the minimum in a sample of size $k$, $F_{1:k}=1-(1-F)^k$. Consequently, the pdf of $U_{m(k)}$ is given
by (see Arnold et al., 1998)
\begin{equation}
 \label{eq:pdf-upper-record}
f_{U_{m(k)}}(u;\theta)
=\frac{k}{(m-1)!}[-\log(1-F(u;\theta))^k]^{m-1}(1-F(u;\theta))^{k-1}f(u;\theta).
\end{equation}

For a comprehensive treatment and for references to the extensive
literature on the topic of $k$-record statistics, one may refer
to the books of Arnold et al. (1998) and Nevzorov (2001).

\section{Main results}
From \eqref{hyb}, the likelihood function of $\theta=(\mu,\sigma)$, under the Rayleigh
distribution in \eqref{eq:cdfre}, is given by
\begin{eqnarray}\nonumber
 \label{eq:pdf-censorII}
L(\theta|\mathbf{x})\propto{\sigma^{-d}}\left[\prod_{i=1}^{d}(x_i-\mu)\right]
\exp\left\{\frac{-1}{2\sigma}\left[\sum_{i=1}^{d}(x_i-\mu)^2+(n-d)(T_0-\mu)^2\right]\right\},
&&
\\ \mu<x_1, \ \  \sigma>0.&&
\end{eqnarray}
Also, from \eqref{eq:pdf-upper-record}, the pdf of $U_{m(k)}$, under the
Rayleigh distribution in \eqref{eq:cdfre}, is
\begin{equation}
\label{eq:pdfrecordexp}
f_{U_{m(k)}}(u|\theta)=\frac{k^m
(u-\mu)^{2m-1}}{\Gamma(m)\sigma^m
2^{m-1}}\exp\left\{-\frac{k(u-\mu)^2}{2\sigma}\right\},~~~~u>\mu.
\end{equation}

In the sequel, we provide the predictive inference, under both
scaled and two parameter Rayleigh distributions.

\subsection{Prediction of a repair time under the scaled Rayleigh distribution}

For the scaled Rayleigh distribution (the case $\mu=0$), we
consider the non-informative prior
\begin{equation}
 \label{eq:pdf1-prior}
\pi_1(\sigma)\propto \frac{1}{\sigma}, \ \ \sigma>0.
\end{equation}

From \eqref{eq:pdf-censorII} and \eqref{eq:pdf1-prior}, the posterior density of $\sigma$ is given
by
\begin{equation}
 \label{eq:pdf-poster1}
\pi^*_1(\sigma|{\mathbf{x}})=\frac{[\delta({\mathbf{x}})]^d}{\Gamma(d)2^d\sigma^{d+1}}\exp\left\{ \frac{-\delta({\mathbf{x}})}{2\sigma}\right\}, \ \  \sigma>0,
\end{equation}
where
\begin{equation}
\label{eq:delta}
\delta({\mathbf{x}})=\sum_{i=1}^{d}x_i^2+(n-d)T_0^2,
\end{equation}
which is the pdf of an inverted
gamma distribution with parameters $d$ and
$\frac{\delta({\mathbf{x}})}{2}$.

From \eqref{eq:pdf-poster1} and \eqref{eq:pdfrecordexp}, the predictive density of $U_{m(k)}$, given
$\mathbf{X}=\mathbf{x}$, is
\begin{eqnarray}
\label{eqpdfr1}
\nonumber
f^*_{U_{m(k)}}(u|\mathbf{x}) &=&\int_0^\infty
f_{U_{m(k)}}(u|\sigma)\pi^*_1(\sigma|\mathbf{x})\textmd{d}\sigma\\ \nonumber
&=&\int_0^\infty \frac{k^m u^{2m-1}[\delta({\mathbf{x}})]^d}{\Gamma(m)\Gamma(d)\sigma^{m+d+1} 2^{d+m-1}}\exp\left\{-\frac{k u^2 +\delta(\mathbf{x})}{2\sigma}\right\}
~\textmd{d}\sigma
\\\nonumber
\\&=&\frac{2}{B(d,m) u} p_k(u;\mathbf{x})^{d}(1-p_k(u;\mathbf{x}))^{m},~~~~u>0,
\end{eqnarray}
where $B(\cdot, \cdot)$ is the complete beta function,
\begin{equation}
\label{eq:p} p_k(y;\mathbf{x})=\frac{\delta(\mathbf{x})}{\delta(\mathbf{x})+k y^2},
\end{equation}
and $\delta(\cdot)$ is given in \eqref{eq:delta}.

From \eqref{eqpdfr1}, it follows that $p_k(U_{m(k)};\mathbf{x})|\mathbf{x}\sim Beta (d,m)$, hence
\begin{eqnarray}
\label{eq:survival}\nonumber
Pr(U_{m(k)}> z|\mathbf{x})&=&Pr(p_k(U_{m(k)};\mathbf{x})<p_k(z;\mathbf{x})|\mathbf{x})
\\ &=&I(d, m, p_{k}(z;\mathbf{x})),\ \ z>0,
\end{eqnarray}
where $I(\gamma_1,\gamma_2,x)$ is the incomplete beta function,
$$I(\gamma_1,\gamma_2,x)=\int_0^x \frac{1}{B(\gamma_1,\gamma_2)}t^{\gamma_1-1}(1-t)^{\gamma_2-1}\mathrm{d}t.$$
\begin{remark}\em
The prior used in \eqref{eq:pdf1-prior} is a special case
of the conjugate inverted gamma prior with pdf
\begin{equation}
 \label{eq:pdf-prior}
\tilde{\pi}(\sigma)\propto
\sigma^{-(a+1)}\exp\left\{-\frac{b}{\sigma}\right\}, \ \ \sigma>0,
\end{equation}
where $a$ and $b$ are positive hyper-parameters. The prior $\tilde{\pi}(\sigma)$ tends to $\pi_1(\sigma)$ as
$(a,b)\rightarrow (0,0)$. Since then the variance of $\sigma$ tends to infinity, we shall call $\pi_1(\sigma)$ the
noninformative conjugate prior for $\sigma$.
\end{remark}

\subsubsection{Interval prediction of a repair time}

A two-sided equi-tailed $100 (1-\alpha)\%$ Bayesian prediction interval (PI) for $U_{m(k)},~m\geq 1$, is
obtained from \eqref{eq:survival}, as the interval $(L(\mathbf{x}),U(\mathbf{x}))$, for which
\begin{eqnarray*}
Pr(U_{m(k)}>L(\mathbf{x})|\mathbf{x})=1-\frac{\alpha}{2} \qquad \mbox{ and } \qquad
Pr(U_{m(k)}>U(\mathbf{x})|\mathbf{x})=\frac{\alpha}{2}.
\end{eqnarray*}
So, the interval $(L(\mathbf{x}),U(\mathbf{x}))$ is
\begin{equation}\label{eq:twotailed-l}
\left(\sqrt{\frac{[1-\beta_{\frac{\alpha}{2}}(d,m)]\delta(\mathbf{x})}{k \beta_{\frac{\alpha}{2}}(d,m)}},\sqrt{\frac{[1-\beta_{1-\frac{\alpha}{2}}(d,m)]\delta(\mathbf{x})}{k \beta_{1-\frac{\alpha}{2}}(d,m)}}\right),
\end{equation}
where $\beta_{\gamma}(n_1,n_2)$ denotes the upper $\gamma^{\rm th}$ quantile of
the beta distribution with parameters $n_1$ and $n_2$, i.e.,
$P(T>\beta_{\gamma}(n_1,n_2))=\gamma$, with $T\sim
Beta(n_1,n_2)$.

The highest posterior density prediction interval (HPD PI) is an interval, the
posterior pdf for every point inside which is
greater than that for every point outside of which. A HPD PI
includes the more probable values of the parameter and excludes
the less probable ones. Since the posterior pdf
$f^*_{U_{m(k)}}(u|\mathbf{x})$ is unimodal and $p_k(u;\mathbf{x})$ is  decreasing in $u$, the HPD PI,
$(w_1,w_2)$, for $U_{m(k)}$ given $\mathbf{X}=\mathbf{x}$, with
coverage probability $1-\alpha$, is the simultaneous solution of 
\begin{equation}\label{hpd1}
\frac{1}{B(d,m)}\int_{p_k(w_2;\mathbf{x})}^{p_k(w_1;\mathbf{x})}t^{d-1}(1-t)^{m-1}\mathrm{d}t=1-\alpha,
\end{equation}
and
\begin{equation}\label{hpd2}
\left(\frac{w_{1}}{w_{2}}\right)^{2m-1}=\left(\frac{\delta(\mathbf{x})+k w^2_{1}}{\delta(\mathbf{x})+k w^2_{2}}\right)^{d+m}.
\end{equation}

\begin{remark}\em
For $m=1$, we have $U_{1(k)}=Y_{1:k}$, that is the lifetime of
the series system, and the prediction interval
\eqref{eq:twotailed-l} is simplified to
 \begin{equation}
\label{eq:twotailed-2}
\left(\sqrt{\frac{\delta(\mathbf{X})}{k}\left\{\left(1-\frac{\alpha}{2}\right)^{-\frac{1}{d}}-1\right\}},
\ \
\sqrt{\frac{\delta(\mathbf{X})}{k}\left\{\left(\frac{\alpha}{2}\right)^{-\frac{1}{d}}-1\right\}}
\right).
\end{equation}
Also, the HPD PI, $(w_1,w_2)$, for $U_{1(k)}$ given
$\mathbf{X}=\mathbf{x}$, with coverage probability $1-\alpha$, is the
simultaneous solution of
\[[p_k(w_1;\mathbf{x})]^{d}-[p_k(w_2;\mathbf{x})]^{d}=1-\alpha\]
and
\[w_1[p_k(w_1;\mathbf{x})]^{d+1}=w_2[p_k(w_2;\mathbf{x})]^{d+1}.\]
\end{remark}

\subsubsection{Point prediction of a repair time}

Using \eqref{eqpdfr1} and under the squared error loss (SEL) function, the Bayes point predictor of $U_{m(k)}$  is
\begin{eqnarray}
\widehat{U}_{m(k)}&=&E(U_{m(k)}|\mathbf{x})\nonumber\\
&=&\int_0^\infty \frac{2}{B(d,m)} p_k(u;\mathbf{x})^{d}(1-p_k(u;\mathbf{x}))^{m}\mathrm{d}u
\nonumber\\&=&\int_0^1 \frac{\sqrt{\frac{\delta(\mathbf{x})}{k}}}{B(d,m)} z^{d-\frac{3}{2}}(1-z)^{m-\frac{1}{2}}\mathrm{d}z\nonumber\\&=&
\frac{B(d-\frac{1}{2},m+\frac{1}{2})}{B(d,m)}\sqrt{\frac{\delta(\mathbf{X})}{k}}.
\end{eqnarray}
Similarly, it can be verified that, under the absolute error loss (AEL) function and zero-one loss function, the Bayes point predictors of $U_{m(k)}$ are  {
$$\hat{U}^*_{m(k)}={\rm Med}(U_{m(k)}|\mathbf{x})=\sqrt{\frac{(1-{\rm{Med}}(d,m))[\delta(\mathbf{X})]}{k{\rm{Med}}(d,m)}},$$}
and
$$\widetilde{U}_{m(k)}={\rm Mod}(U_{m(k)}|\mathbf{x})=\sqrt{\frac{(2m-1)[\delta(\mathbf{X})]}{k(2d+1)}},$$
respectively, where {Med$(d,m)$} denotes the median of Beta distribution with parameters $d$ and $m$.

\subsection{Results for the two parameter Rayleigh distribution}
To facilitate the Bayesian approach under the two parameter Rayleigh distribution,
we assume independent prior distributions for the model
parameters, that is
\begin{equation}\label{prior2}
\pi(\mu,\sigma)\propto \pi_1(\sigma)\pi_2(\mu),
\end{equation}
where $\pi_1(\sigma)$ is the non-informative conjugate prior in
\eqref{eq:pdf1-prior} and  $\pi_2(\mu)$ is a normal density with mean
$\xi$ and variance $1/2\tau$, that is
\begin{equation}
 \label{eq:pdf-prior-mu}
\pi_2(\mu)\propto \exp\left\{-\tau(\mu-\xi)^2\right\}, \ \
\mu\in\mathbb{R}, \ \ \xi\in\mathbb{R},\ \ \tau>0.
\end{equation}
Therefore, the joint prior density is
\begin{equation}
 \label{eq:pdf-prior-mu-sigma}
\pi(\mu,\sigma)\propto \frac{1}{\sigma}
\exp\left\{-\tau(\mu-\xi)^2\right\},
 \ \ \sigma>0, \ \ \mu\in\mathbb{R}, \ \ \xi\in\mathbb{R},\ \
 \tau>0,
\end{equation}
and the posterior density function of $\mu$ and
$\sigma$, given ${\bf {X=x}}$, is obtained as
\begin{eqnarray}\nonumber
 \label{eq:pdf-poster}
\pi^*(\mu,
\sigma|{\mathbf{x}})=\frac{A_1(\mathbf{x})\exp\left\{-\tau(\mu-\xi)^2\right\}}{\Gamma(d)2^{d}\sigma^{d+1}}
\left[\prod_{i=1}^d(x_i-\mu)\right]\exp\left\{
\frac{-1}{2\sigma}[\delta^*(\mu|{\mathbf{x}})]\right\},&& \\
 \sigma>0, \ \ \mu<x_1,&&
\end{eqnarray}
where
\begin{equation}
 \label{eq:delta-star}
 \delta^*(\mu|{\mathbf{x}})=\sum_{i=1}^{d}(x_i-\mu)^2+(n-d)(T_0-\mu)^2
\end{equation}
and
\begin{equation}
 \label{eq:A1}
A_1(\mathbf{x})^{-1} =\int_{-\infty}^{x_{1}}
\frac{\exp\left\{-\tau(\mu-\xi)^2\right\}\prod_{i=1}^d(x_i-\mu)}{[\delta^*(\mu|{\mathbf{x}})]^{d}}
\textrm{d}\mu.
\end{equation}
The predictive posterior density of $U_{m(k)}$, given
$\mathbf{x}$, is obtained as follows
\begin{eqnarray}
\hspace{-1.5 cm} \nonumber
h^*_{U_{m(k)}}(u|\mathbf{x}) &=&\int_{-\infty}^{\min(u,x_1)}\int_0^\infty
f_{U_{m(k)}}(u|\mu,\sigma)\pi^*(\mu,\sigma|\mathbf{x})\textmd{d}\sigma\textmd{d}\mu\\
\nonumber \\
&=&\frac{2A_1(\mathbf{x})}{B(d,m)} \int_{-\infty}^{\min(u,x_1)}
\frac{g(t,u,m;\mathbf{x},{\boldsymbol\eta})}{u-t}\textrm{d}t,\quad
u\in\mathbb{R},\label{eqpdfr22}
\end{eqnarray}
where ${\boldsymbol\eta}=(\tau,\xi,k,d)$ and
\begin{equation}\label{gf}
g(t,u,m;\mathbf{x},{\boldsymbol\eta})=\frac{k^m\exp\left\{-\tau(t-\xi)^2\right\}
(u-t)^{2m}
\prod_{i=1}^d(x_i-t)}{[k(u-t)^2+\delta^*(t|{\mathbf{x}})]^{d+m}}.
\end{equation}
%The proof of \eqref{eqpdfr22} is given in the Appendix.

The predictive posterior survival function of $U_{m(k)}$ is, for $z\geq x_1$,
\begin{eqnarray}
\label{eq:survival2}
\bar{H}^*_{U_{m(k)}}(z|\mathbf{x})&=&
\int_{-\infty}^{x_1}\int_0^\infty\int_z^\infty f_{U_{m(k)}}(u|\mu,\sigma)\pi_2^*(\mu,\sigma)
\textrm{d}u\textmd{d}\sigma \textrm{d}\mu\nonumber\\
&=&A_1(\mathbf{x})\sum_{j=0}^{m-1}{\frac{\Gamma(d+j)}{\Gamma(d)j!}}\int_{-\infty}^{x_1} g(t,z,j;\mathbf{x},{\boldsymbol\eta})\textrm{d}t,
\end{eqnarray}
where $g(t,z,j;\mathbf{x},{\boldsymbol\eta})$ is given in \eqref{gf}; and for $z<x_1$,
\begin{eqnarray}
\label{eq:survival3}\bar{H}^*_{U_{m(k)}}(z|\mathbf{x})&=&A_1(\mathbf{x})\left[\sum_{j=0}^{m-1}{\frac{\Gamma(d+j)}{\Gamma(d)j!}}
\int_{-\infty}^{z}g(t,z,j;\mathbf{x},{\boldsymbol\eta})\;{\rm d}t\right.\nonumber\\
&&\left.+\int_{z}^{x_1}g(t,t,0;\mathbf{x},{\boldsymbol\eta})\;\textrm{d}t\right],
\end{eqnarray}
wherein $0^0$ is defined to be 1 and $g$ is given in \eqref{gf}. 

%The proofs of \eqref{eq:survival2} and \eqref{eq:survival3} are given in the Appendix.

\subsubsection{Interval prediction of a repair time}

The equi-tailed $100(1-\alpha)\%$ Bayesian prediction interval for $U_{m(k)}$ can be obtained numerically, using \eqref{eq:survival2} and \eqref{eq:survival3}.

{Since $h^*_{U_{m(k)}}(u|\mathbf{x})$ is unimodal,} then the HPD PI, $(w_1,w_2)$, for $U_{m(k)}$, with coverage probability $1-\alpha$,  satisfies 
$$\int_{w_1}^{w_2} h^*_{U_{m(k)}}(u|\mathbf{x})du=1-\alpha$$
and
$$h^*_{U_{m(k)}}(w_1|\mathbf{x})=h^*_{U_{m(k)}}(w_2|\mathbf{x}),$$
where $h^*_{U_{m(k)}}(u|\mathbf{x})$ is given in \eqref{eqpdfr22}.

\subsubsection{Point prediction of a repair time}

Under the SEL function, the point predictor of $U_{m(k)}$ is given by
\begin{eqnarray}\label{selpp2}
\widehat{U}_{m(k)}&=&E(U_{m(k)}|\mathbf{x})\nonumber\\
&=& A_1(\mathbf{x})\left[\frac{\Gamma(m+1/2)\Gamma(d-1/2)}{\Gamma(m)\Gamma(d) \sqrt{k}}\int_{-\infty}^{x_1}\sqrt{\delta^*(t|\mathbf{x})}g(t,t,0;\mathbf{x},{\boldsymbol\eta})\;{\rm d}t\right.\nonumber\\
&&\left.+\int_{-\infty}^{x_1}t\;g(t,t,0;\mathbf{x},{\boldsymbol\eta})\;\textrm{d}t\right],
\end{eqnarray}
where $g$ is given in \eqref{gf}. 

%The proof of \eqref{selpp2} is given in the Appendix.

For the AEL function, ${U^*}_{m(k)}$ can be obtained as 
$${U^*}_{m(k)}={{\bar{H^*}}^{-1}}_{U_{m(k)}}\left(\frac{1}{2}\right).$$
Also, the point predictor $\tilde{U}_{m(k)}$, under the zero-one loss function is the unique mode of the pdf $h^*_{U_{m(k)}}(w|\mathbf{x})$
in \eqref{eqpdfr22}.

\subsection{Model checking}

The assumption of independence of $\pi_1(\sigma)$ and $\pi_2(\mu)$
in \eqref{prior2} may affect on the performance of the predictors.
Furthermore, as $\tau$ tends to zero, the prior $\pi_2(\mu)$ tends
to a noninformative prior, while the mean square error of prediction,
$$\E(\widehat{U}_{m(k)}-{U}_{m(k)})^2,$$
increases. Hence, choosing a suitable value for $\tau$ is an
important issue. It should be kept in mind that $\widehat{U}_{m(k)}$ minimizes 
$$\E[(g(\mathbf{x})-{U}_{m(k)})^2|\mathbf{x}),$$
over all functions $g$, while $\E(\widehat{U}_{m(k)}-{U}_{m(k)})^2$ can be greater or less than $\E({U'}_{m(k)}-{U}_{m(k)})^2=2\var({U}_{m(k)})$, in which ${U}_{m(k)}$ and ${U'}_{m(k)}$ are iid random variables from the distribution 
\eqref{eq:pdfrecordexp}. Hence,for the prior distribution in \eqref{prior2} to be 
sufficiently low informative, and for the mean square error of prediction not to be very large, one may choose $\tau$ such that 
$$\frac{\E(\widehat{U}_{m(k)}-{U}_{m(k)})^2}{\E({U'}_{m(k)}-{U}_{m(k)})^2}$$
takes its largest value less than or equal to 1. 

To check the suitableness of the prior distribution
in \eqref{prior2} and to choose a suitable value for $\tau$, we
perform a simulation study as follows. This study is based on the
general method for model checking, described by Gelman et al. (2004).

The algorithm for model checking is as follows:

\noindent{\bf Algorithm:}

\begin{description}
\item (i) Generate  $x_1,\cdots,x_n$ independently from Rayleigh distribution and extract the hybrid censored sample $(x_1,\cdots,x_d)$ from  $x_1,\cdots,x_n$.  
\item (ii) Generate $u_{m(k)}$ and $u'_{m(k)}$ independently from \eqref{eq:pdfrecordexp}.
\item (iii) Predict $u^*_{m(k)}$ using \eqref{selpp2}.
\item (iv) Replicate (i)-(iii), $N=10000$ times, independently, to obtain
samples $u_{m(k)}^{(1)},\ldots,u_{m(k)}^{(N)}$ and ${u'}^{(1)}_{m(k)},\ldots,{u'}^{(N)}_{m(k)}$ as well as predicted values
${u^*}^{(1)}_{m(k)},\ldots,{u^*}^{(N)}_{m(k)}$.
\item (v) Compute
$$SS_1=\frac{1}{N}\sum_{i=1}^{N}(u_{m(k)}^{(i)}-{u^*}^{(i)}_{m(k)})^2$$
and
$$SS_2=\frac{1}{N}\sum_{i=1}^{N}(u_{m(k)}^{(i)}-{u'}^{(i)}_{m(k)})^2.$$
\item (vi) Set $\xi=0$ and choose $l^*$ from the set $\{-2,-1,\ldots,2\}$, such that, for  $\tau^*=0.5\times 10^{-l^*}$, $D_1=SS_1/SS_2$ takes its largest value less than or equal to 1.
\item (vii) For $\xi=0$ and $\tau=\tau^*$, compute 
$$D_2=\frac{\frac{1}{N}\sum_{i=1}^{N}(u_{m(k)}^{(i)}-\mu)}{\frac{1}{N}\sum_{i=1}^{N}({u^*}^{(i)}_{m(k)}-\mu)}$$
and 
$$D_3=\frac{\frac{1}{N}\sum_{i=1}^{N}(u_{m(k)}^{(i)})^2-\left(\frac{1}{N}\sum_{i=1}^{N}u_{m(k)}^{(i)}\right)^2}{\frac{1}{N}\sum_{i=1}^{N}({u^*}^{(i)}_{m(k)})^2-\left(\frac{1}{N}\sum_{i=1}^{N}{u^*}^{(i)}_{m(k)}\right)^2}.$$
\end{description}

\begin{figure}
\centering 
\includegraphics[scale=0.29]{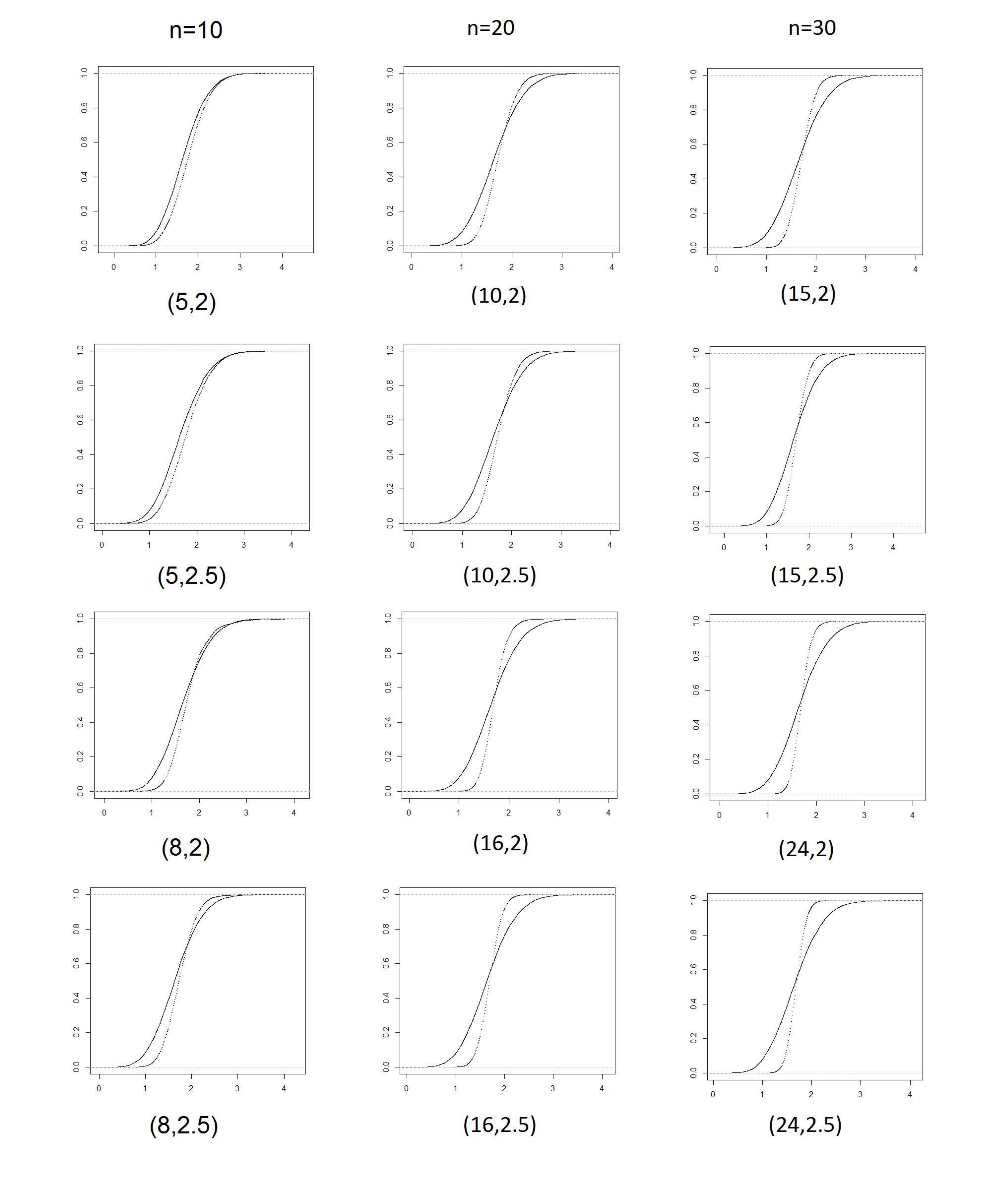}
\caption{Empirical cdfs of the simulated values, $u_{m(k)}^{(1)},\ldots,u_{m(k)}^{(N)}$ (solid line), and the predicted values, ${u^*}^{(1)}_{m(k)},\ldots,{u^*}^{(N)}_{m(k)}$ (dashed line), for  $m=3$, $k=2$, $\mu=0$, $\sigma=1$, $\xi=0$, $\tau=\tau^*$ and different values of $n$ and $(r,T)$. }\label{ecdf}
\end{figure}

Table \ref{mc} presents the values of $\tau^*$, $D_1$, $D_2$ and $D_3$, for $m=3$, $k=2$, $\mu=0$, $\sigma=1$, $\xi=0$ and $\tau=\tau^*$.
Also, empirical cdfs of the simulated values, $u_{m(k)}^{(1)},\ldots,u_{m(k)}^{(N)}$, and the predicted values, ${u^*}^{(1)}_{m(k)},\ldots,{u^*}^{(N)}_{m(k)}$, are shown in Figure \ref{ecdf}, for  $m=3$, $k=2$, $\mu=0$, $\sigma=1$, $\xi=0$, $\tau=\tau^*$ and different values of $n$ and $(r,T)$. 
As one can see from Table \ref{mc} and Figure \ref{ecdf}, for a sufficiently low informative prior, the predicted values have smaller variance and mean square error of prediction than the simulated values. The ignorable bias of the predicted values decreases as $n$ and/or $r$ get large. Hence, the prior \eqref{prior2} results in efficient predictors. 

\begin{table}
  \centering
  \caption{The values of $\tau^*$, $D_1$, $D_2$ and $D_3$, for $m=3$, $k=2$, $\mu=0$, $\sigma=1$, $\xi=0$ and $\tau=\tau^*$.}\label{mc}
  \begin{tabular}{c c c c c c}
    \hline
    $n$ & $(r,T)$ & $\tau^*$ & $D_1$& $D_2$ & $D_3$  \\
    10 & (5,2)  &  0.5 & 0.9426 &  0.9330 & 1.2006 \\
       & (5,2.5)  & 0.5 & 0.9417 & 0.9348 & 1.2513 \\
       & (8,2)    & 0.005 & 0.8605  & 0.9441  & 1.5197\\
       & (8,2.5)  & 0.005 & 0.7295  & 0.9508  & 2.2776 \\
    20 & (10,2) & 0.05 & 0.7160  & 0.9577  & 2.5605\\
       & (10,2.5)   & 0.005 & 0.7143  & 0.9561  & 2.4767 \\
       & (16,2)   & 0.005 & 0.6188  &  0.9750 & 4.3192 \\
        & (16,2.5) & 0.005 & 0.6125  &  0.9806 & 5.1078 \\
    30 & (15,2) & 0.05 & 0.6302  & 0.9740  & 4.2603 \\
       & (15,2.5)   & 0.005 & 0.6262  & 0.9758 & 4.1887 \\
       & (24,2)   & 0.005 & 0.5595  & 0.9852  &  6.9719 \\
        & (24,2.5) & 0.005 &  0.5714 &  0.9858 &  7.8098 \\
    \hline
  \end{tabular}
\end{table}

\section{Illustrative example: groove ball bearings data set}

In this section, we illustrate the proposed procedures in the previous section, using the following real data example.

\begin{exam}\label{exam1}
{\rm Consider the endurance test of deep groove ball bearings,
discussed by Leiblein and Zelen (1956), which includes the number
of revolutions (in hundreds of millions) to failure, for each of
$n = 23$ ball bearings, as follows:}
\end{exam}
 \begin{verbatim}
0.1788 0.2892 0.3300 0.4152 0.4212 0.4560 0.4880 0.5184 0.5196 0.5412
0.5556 0.6780 0.6864 0.6864 0.6888 0.8412 0.9312 0.9864 1.0512 1.0584
1.2792 1.2804  1.7340.
\end{verbatim}
The adequacy of the fitness
of the two parameter Rayleigh distribution with
$\hat{\mu}=0.1788$ and $\hat{\sigma}=0.2149$ to the data is tested using the
Kolmogorov-Smirnov (K-S) test. The value of K-S test
statistic is obtained as $D=0.1982$ with a
corresponding $p$-value=0.3269. Hence, the two parameter Rayleigh
distribution fits the data quite well.

Consider the following two sampling schemes:
\begin{itemize}
\item  Scheme 1:~~~~ $r=20$~~~~~and~~~~$T=1.25.$

We have, $d=r=20$, $T_0=\min\{x_{r:n},T\}=1.0584$, and the hybrid censored sample is
\begin{verbatim}
0.1788 0.2892 0.3300 0.4152 0.4212 0.4560 0.4880 0.5184 0.5196 0.5412
0.5556 0.6780 0.6864 0.6864 0.6888 0.8412 0.9312 0.9864 1.0512 1.0584.
 \end{verbatim}\vspace{-1cm}
\item  Scheme 2:~~~~ $r=20$~~~~~and~~~~$T=1.$

We have, in this case, $d=18$, $T_0=1$ and the hybrid censored sample is
\begin{verbatim}
0.1788 0.2892 0.3300 0.4152 0.4212 0.4560 0.4880 0.5184 0.5196
0.5412 0.5556 0.6780 0.6864 0.6864 0.6888 0.8412 0.9312 0.9864.
\end{verbatim}\vspace{-.1cm}
\end{itemize}

\begin{figure}[!hbtp]
\includegraphics[scale=0.6]{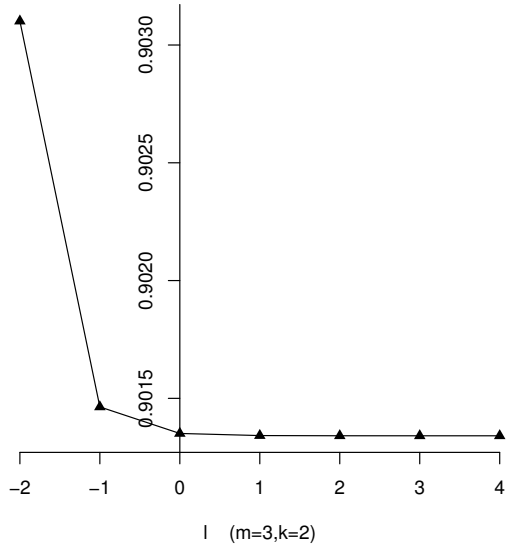}\includegraphics[scale=0.6]{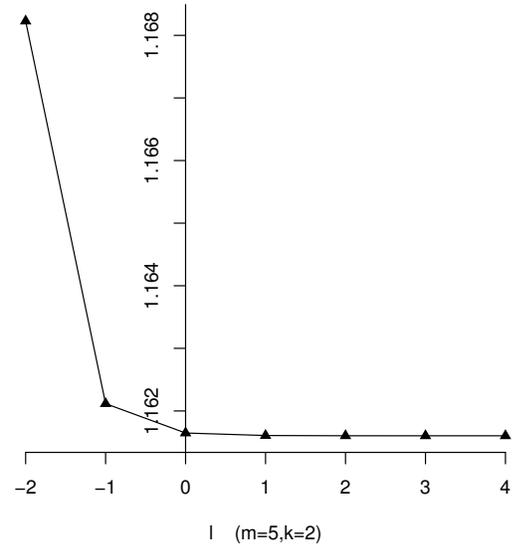}\\
\includegraphics[scale=0.6]{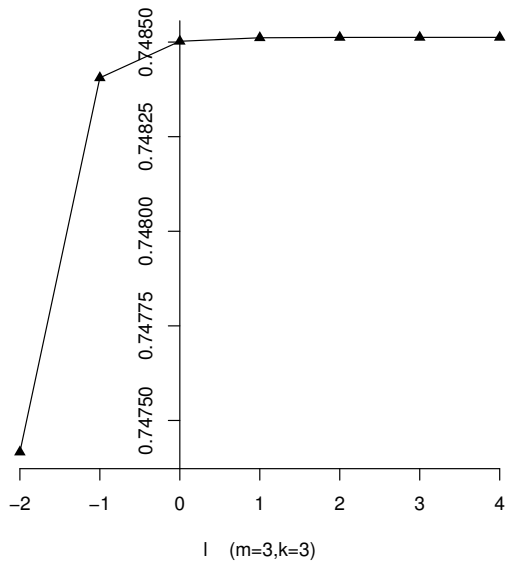}\includegraphics[scale=0.6]{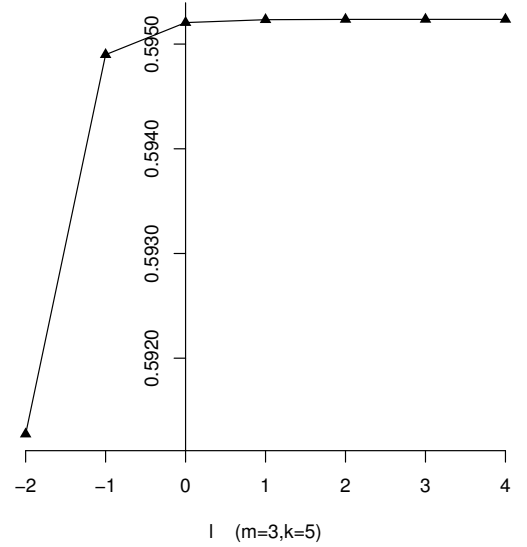}\\
%\centerline{\psfig{figure=sens1.eps,width=6cm}\psfig{figure=sens2.eps,width=6cm}}
%\centerline{\psfig{figure=sens3.eps,width=6cm}\psfig{figure=sens4.eps,width=6cm}}
\caption{\scriptsize Sensitivity analysis for different values of $m$ and $k$.
Plot of the posterior mean versus values of $l$ for $\tau=0.5\times10^{-l}$.}\label{sens}
\end{figure}

To obtain stable Bayesian predictors, we need to
choose suitable values for the hyper-parameters $\xi$ and $\tau$, which ensures that the prior distribution is sufficiently low-informative. This is performed through a sensitivity analysis in the next subsection.

\subsection{Sensitivity analysis}

In order to check the effect of hyper-parameters $\xi$ and $\tau$
on Bayesian predictors, we consider a pilot run on the prediction
procedure, for various values of $\xi$ and $\tau$. As
$\tau\rightarrow 0$, the variance of the normal prior for the
parameter $\mu$ tends to infinity, which is the non-informative
case. As the variance of the normal prior increases, the effect of the hyper-parameter $\xi$ on the
predictors decreases. Therefore, many authors perform the sensitivity analysis on the hyper-parameter $\tau$,
assuming $\xi$ to be fixed as $\xi=0$. So, we take $\xi=0$ and focus on a
sequence of values of $\tau$, tending to 0. Figure \ref{sens} shows
the plot of $\widehat{U}_{m(k)}$ versus values of $l$, for
$\tau=0.5\times10^{-l}$, and different values of $m$ and $k$, when
$\xi=0$. Figure \ref{sens} shows
that for $\tau\geq 0.5$, the values of
$\widehat{U}_{m(k)}$ tend to be very close to each
other. Therefore, $\tau=0.5$ and $\xi=0$ are suitable hyper-parameters for obtaining stable
Bayesian predictors.

\subsection{Bayesian prediction} {
\begin{table}\scriptsize
\center{\caption{\label{tab1}\textmd{95\% Bayesian prediction intervals and point
predictors for $U_{m(k)}$, based on the hybrid censored sample in Example \ref{exam1}}.}}
{\centering Scheme 1}\\\vspace{0.1cm}
\begin{tabular}{c c c c c c c}
\hline
       &$m$&          Equi-tailed PI [width]  &  HPD PI  [width] & $\widehat{U}_{m(k)}$ & $\hat{U}^*_{m(k)}$ & $\tilde{U}_{m(k)}$\\
       \hline
       &1&     (0.1575,1.5235)    [1.3660]   &(0.1004,1.4219)   [1.3215] &0.7111&0.6663&0.5692 \\
       &2&     (0.4070,1.8904)    [1.4834]   &(0.3468,1.7975)   [1.4507] &1.0402&1.0027&0.9302 \\
 $k=1$ &3&     (0.6175,2.1717)    [1.5542]   &(0.5567,2.0811)   [1.5244] &1.2870&1.2498&1.1796 \\
       &4&     (0.7965,2.4113)    [1.6148]   &(0.7337,2.3184)   [1.5847] &1.4928&1.4542&1.3825 \\
\hline
\hline
       &1&     (0.1105,1.0894)    [0.9789]   &(0.0736,1.0329)   [0.9593] &0.5183&0.4895&0.4260 \\
       &2&     (0.2966,1.3448)    [1.0482]   &(0.2580,1.2885)   [1.0305] &0.7510&0.7268&0.6798 \\
 $k=2$ &3&     (0.4495,1.5408)    [1.0913]   &(0.4099,1.4839)   [1.0740] &0.9255&0.9014&0.8556 \\
       &4&     (0.5790,1.7078)    [1.1288]   &(0.5377,1.6485)   [1.1108] &1.0710&1.0459&0.9989 \\
\hline
\hline
       &1&     (0.0851,0.8988)    [0.8137]   &(0.0584,0.8616)   [0.8032] &0.4329&0.4114&0.3632 \\
       &2&     (0.2445,1.1051)    [0.8606]   &(0.2160,1.0656)   [0.8496] &0.6229&0.6048&0.5693 \\
 $k=3$ &3&     (0.3724,1.2633)    [0.8909]   &(0.3424,1.2219)   [0.8795] &0.7654&0.7472&0.7124 \\
       &4&     (0.4801,1.3982)    [0.9181]   &(0.4485,1.3542)   [0.9057] &0.8842&0.8650&0.8291 \\
\hline

\end{tabular}\\\vspace{0.3cm}
{\centering Scheme 2}\\\vspace{0.1cm}
\begin{tabular}{c c c c c c c }
\hline
       &$m$&         Equi-tailed PI [width]  &  HPD PI  [width] & $\widehat{U}_{m(k)}$ & $\hat{U}^*_{m(k)}$ & $\tilde{U}_{m(k)}$\\
       \hline
       &1&     (0.1487,1.5871)   [1.4384]    &(0.0880,1.4792)  [1.3912]  &0.7291&0.6811&0.5785 \\
       &2&     (0.4101,1.9789)   [1.5688]    &(0.3444,1.8754)  [1.5310]  &1.0737&1.0320&0.9526 \\
 $k=1$ &3&     (0.6290,2.2807)   [1.6517]    &(0.5612,2.1767)  [1.6155]  &1.3322&1.2898&1.2112 \\
       &4&     (0.8143,2.5383)   [1.7240]    &(0.7435,2.4299)  [1.6864]  &1.5476&1.5031&1.4215 \\
\hline
\hline
       &1&     (0.0980,1.1290)   [1.0310]    &(0.0592,1.0700)  [1.0108]  &0.5272&0.4967&0.4302 \\
       &2&     (0.2947,1.4011)   [1.1064]    &(0.2527,1.3390)  [1.0863]  &0.7709&0.7442&0.6931 \\
 $k=2$ &3&     (0.4543,1.6108)   [1.1565]    &(0.4102,1.5461)  [1.1359]  &0.9536&0.9263&0.8754 \\
       &4&     (0.5888,1.7901)   [1.2013]    &(0.5422,1.7212)  [1.1790]  &1.1060&1.0771&1.0238 \\
\hline
\hline
       &1&     (0.0706,0.9282)   [0.8576]    &(0.0427,0.8900)  [0.8473]  &0.4377&0.4153&0.3653 \\
       &2&     (0.2398,1.1474)   [0.9076]    &(0.2090,1.1045)  [0.8955]  &0.6367&0.6169&0.5786 \\
 $k=3$ &3&     (0.3737,1.3164)   [0.9427]    &(0.3405,1.2697)  [0.9292]  &0.7860&0.7655&0.7269 \\
       &4&     (0.4860,1.4609)   [0.9749]    &(0.4504,1.4101)  [0.9597]  &0.9103&0.8884&0.8479 \\
\hline

\end{tabular}
\end{table}
}

Here, for the sake of comparison of the performance of the predictors, we take $k=1,2,3$ and $m=1,\ldots,4$.

\begin{figure}[!hbtp]
\includegraphics[scale=0.4]{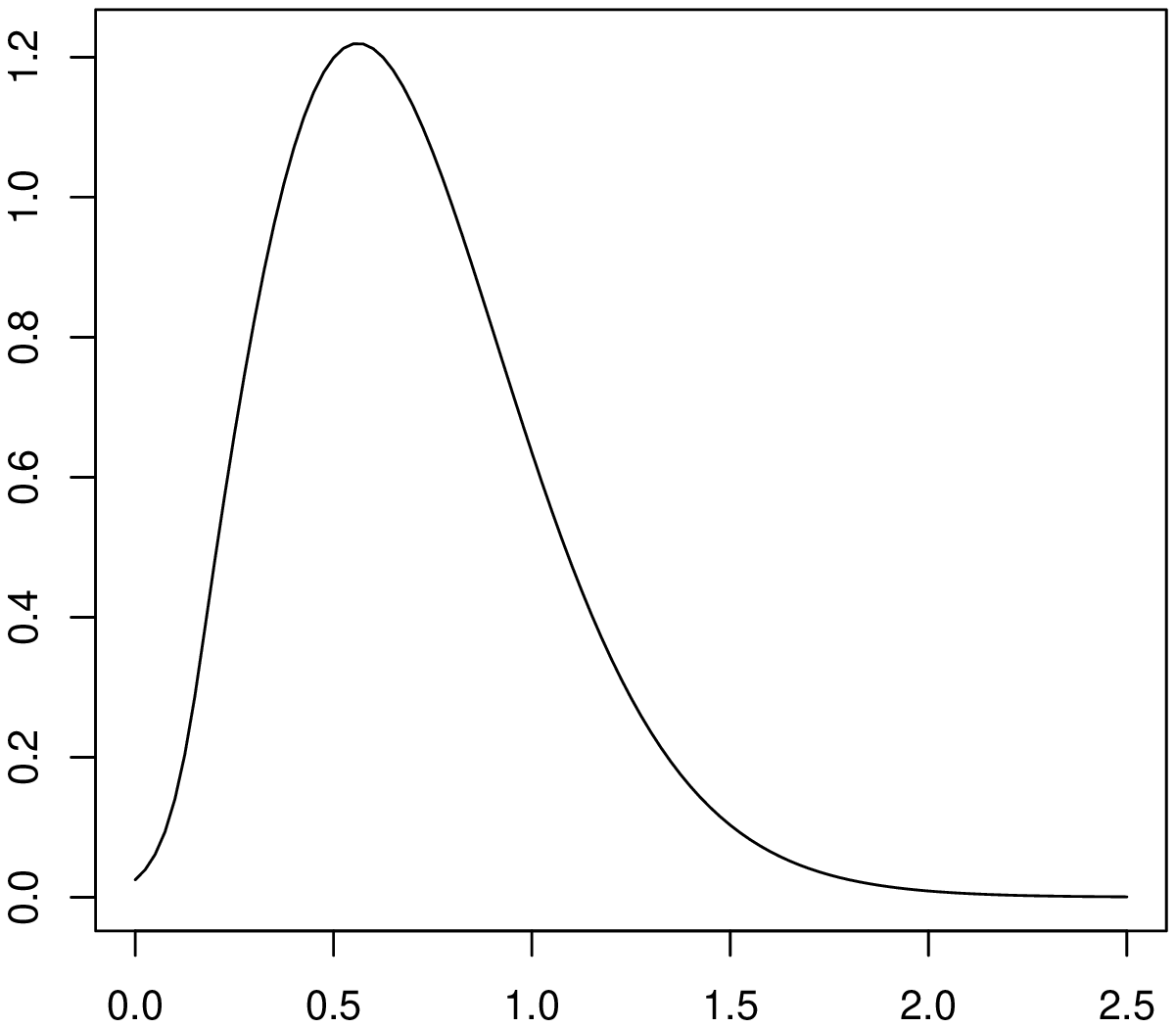}\includegraphics[scale=0.4]{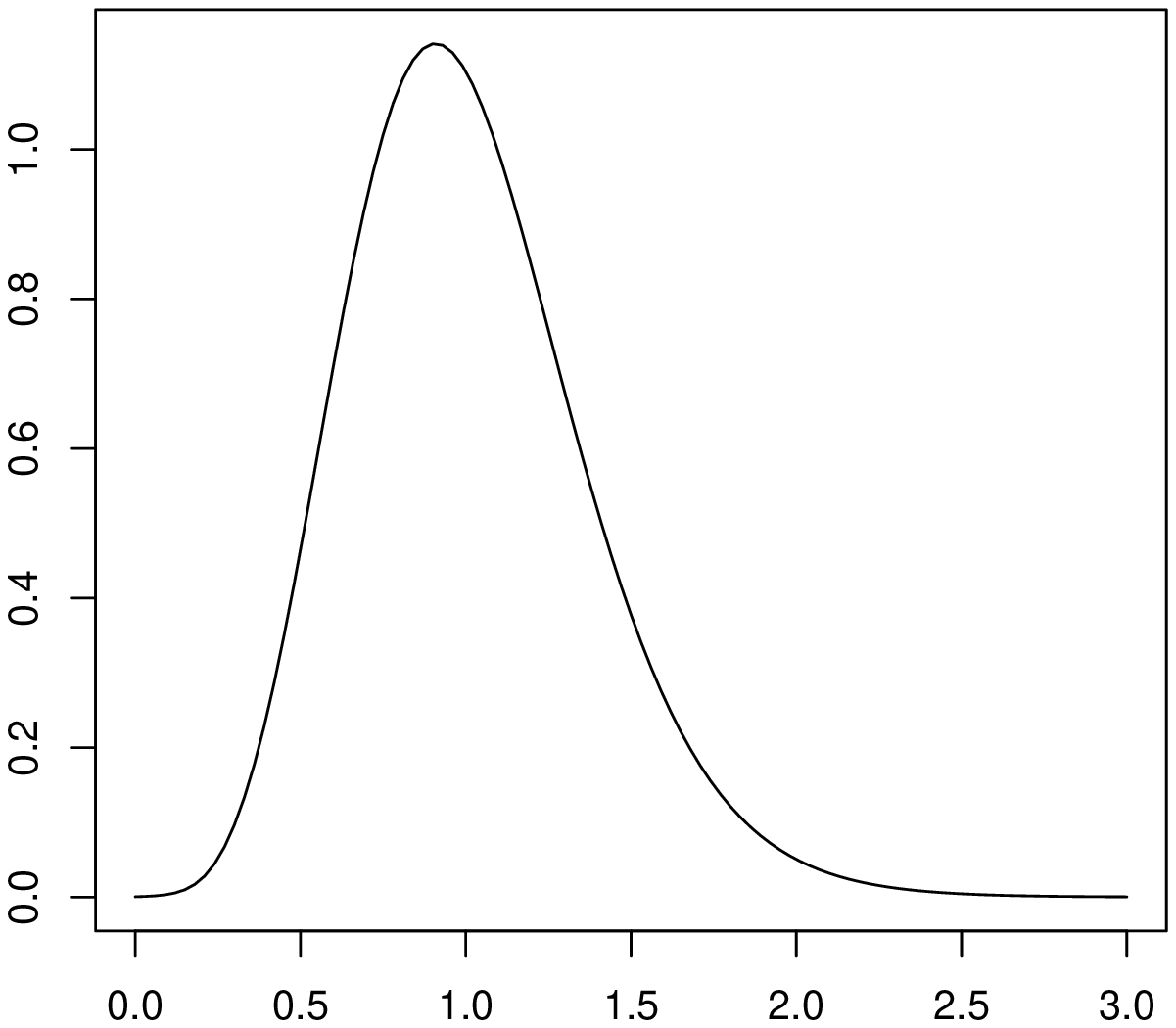}\includegraphics[scale=0.4]{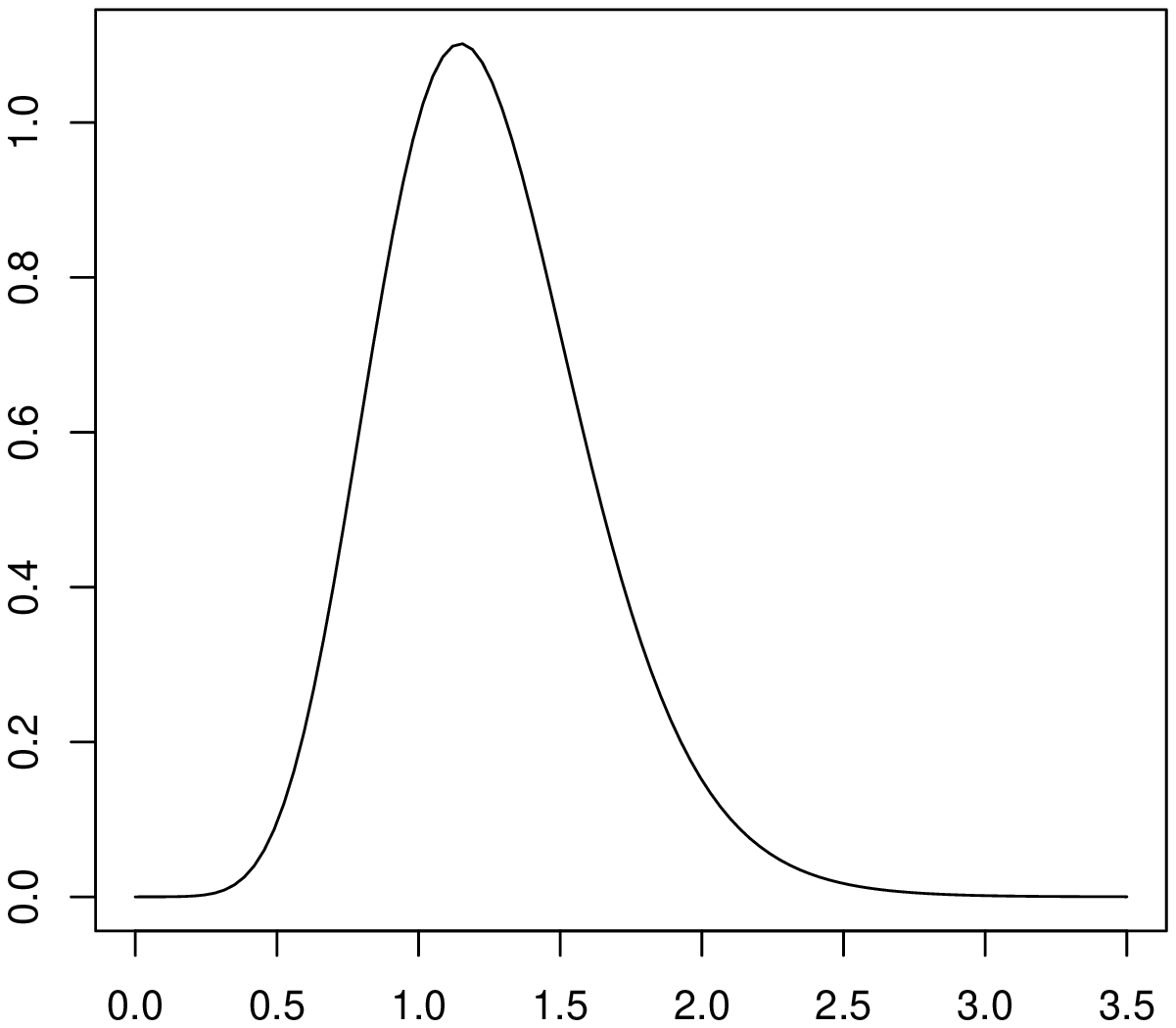}\\
\includegraphics[scale=0.4]{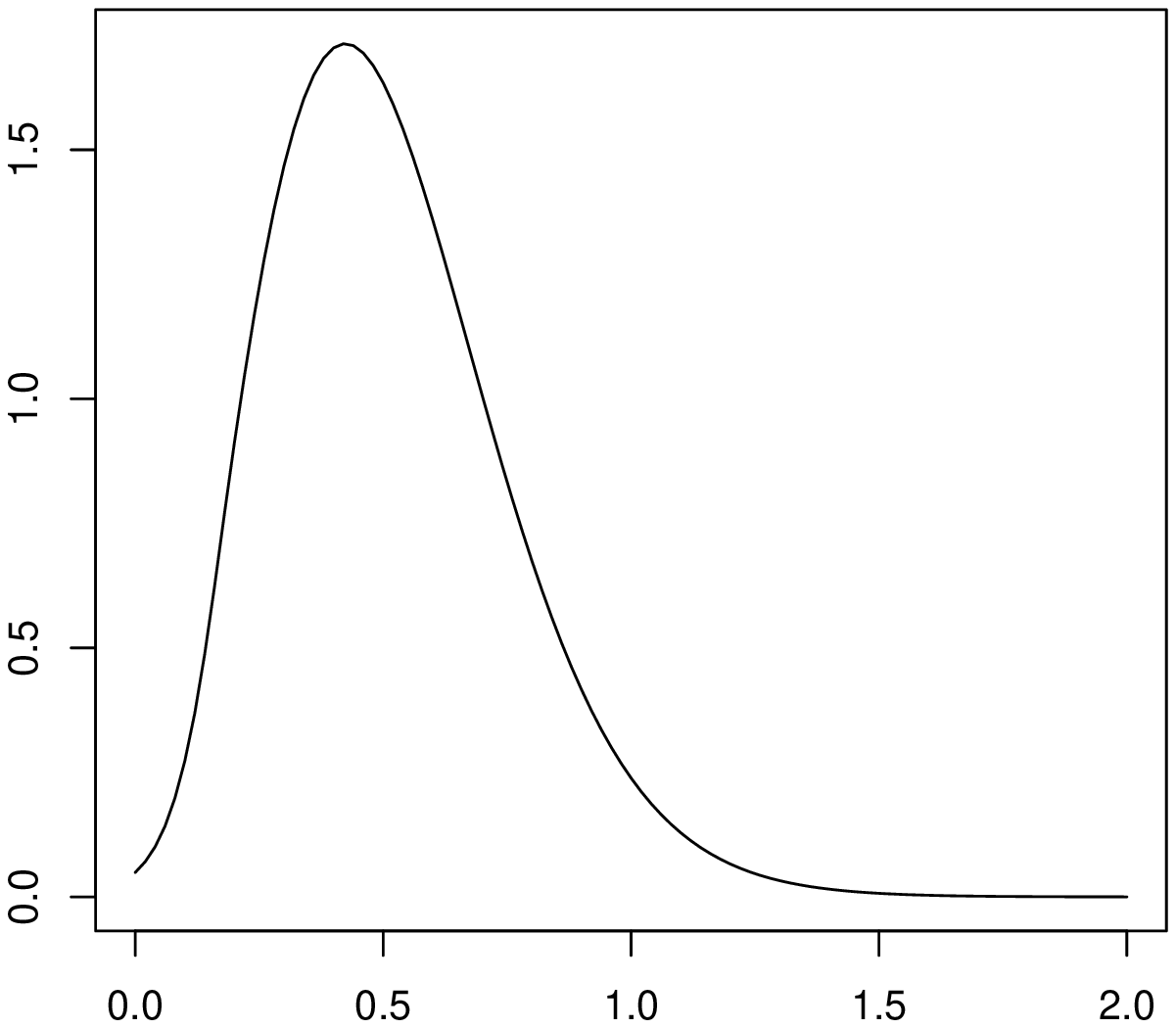}\includegraphics[scale=0.4]{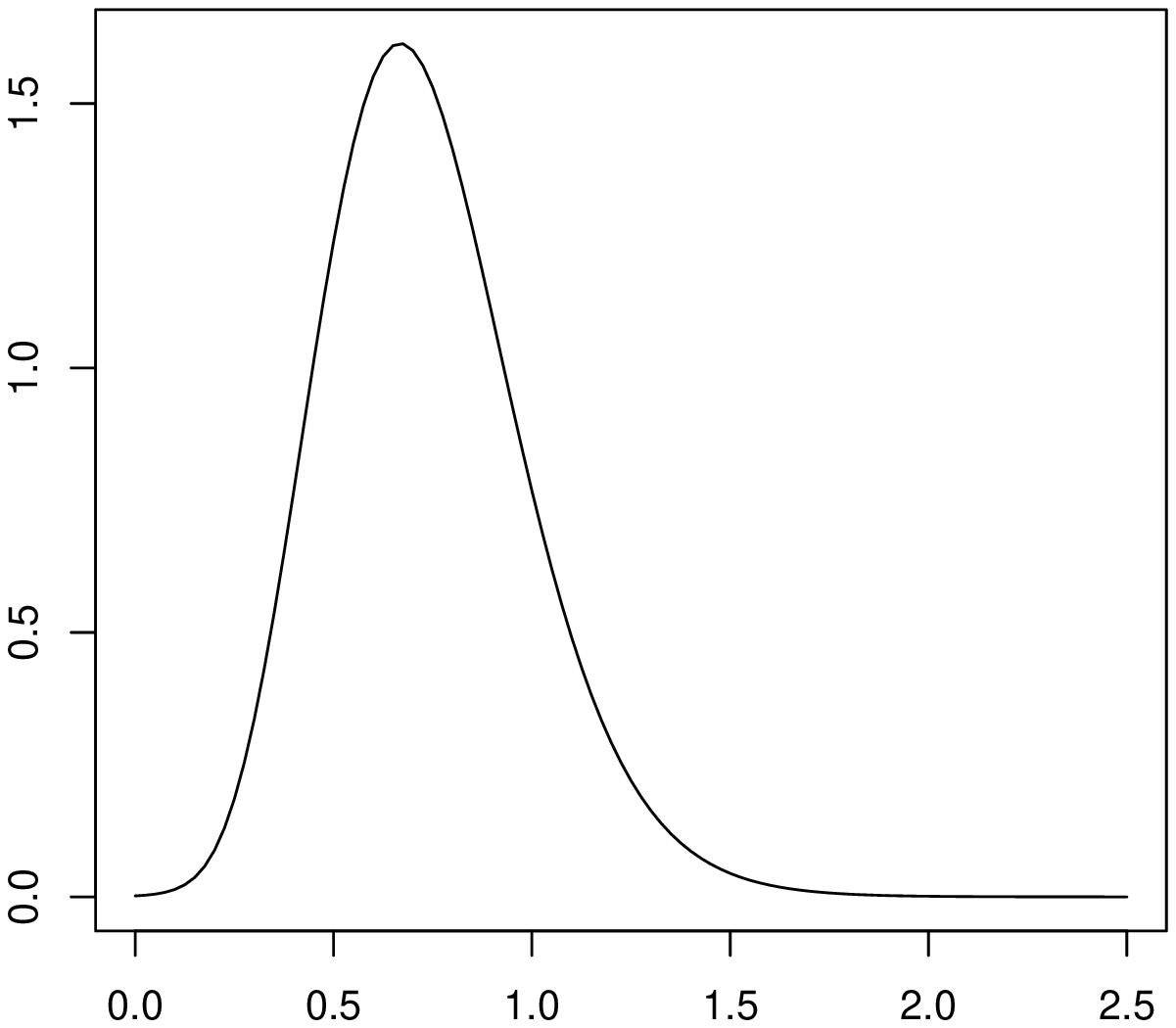}\includegraphics[scale=0.4]{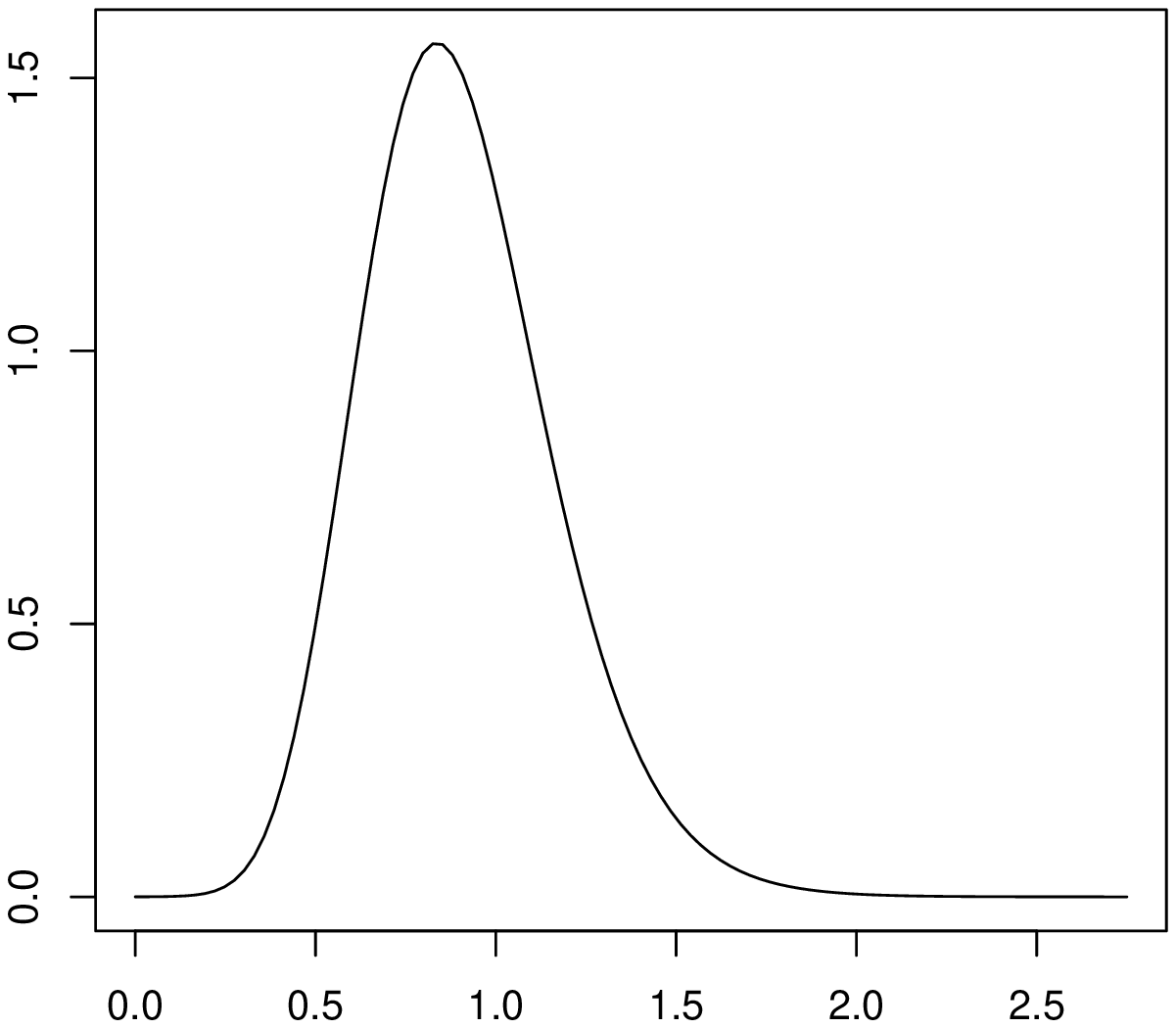}\\
\includegraphics[scale=0.4]{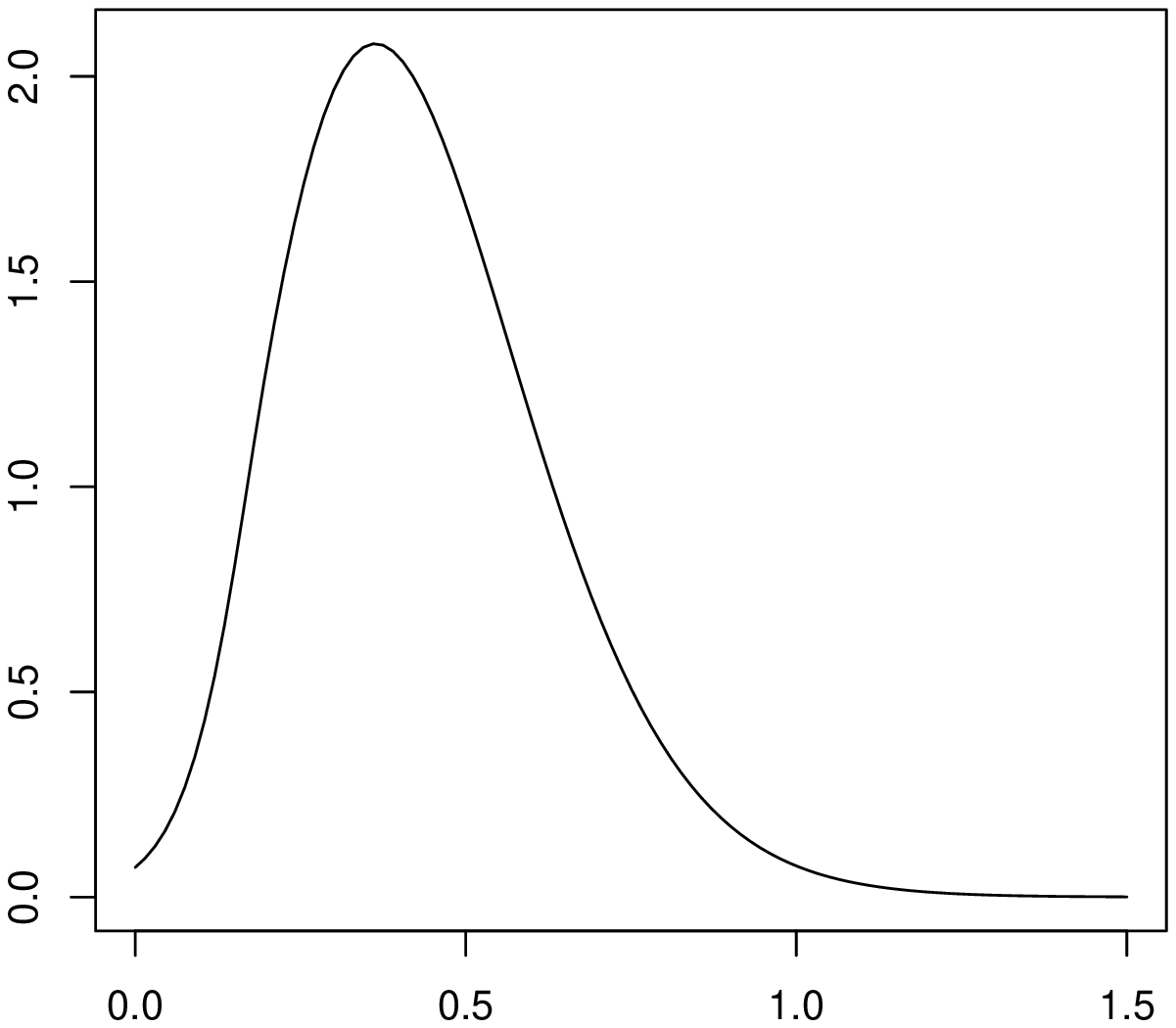}\includegraphics[scale=0.4]{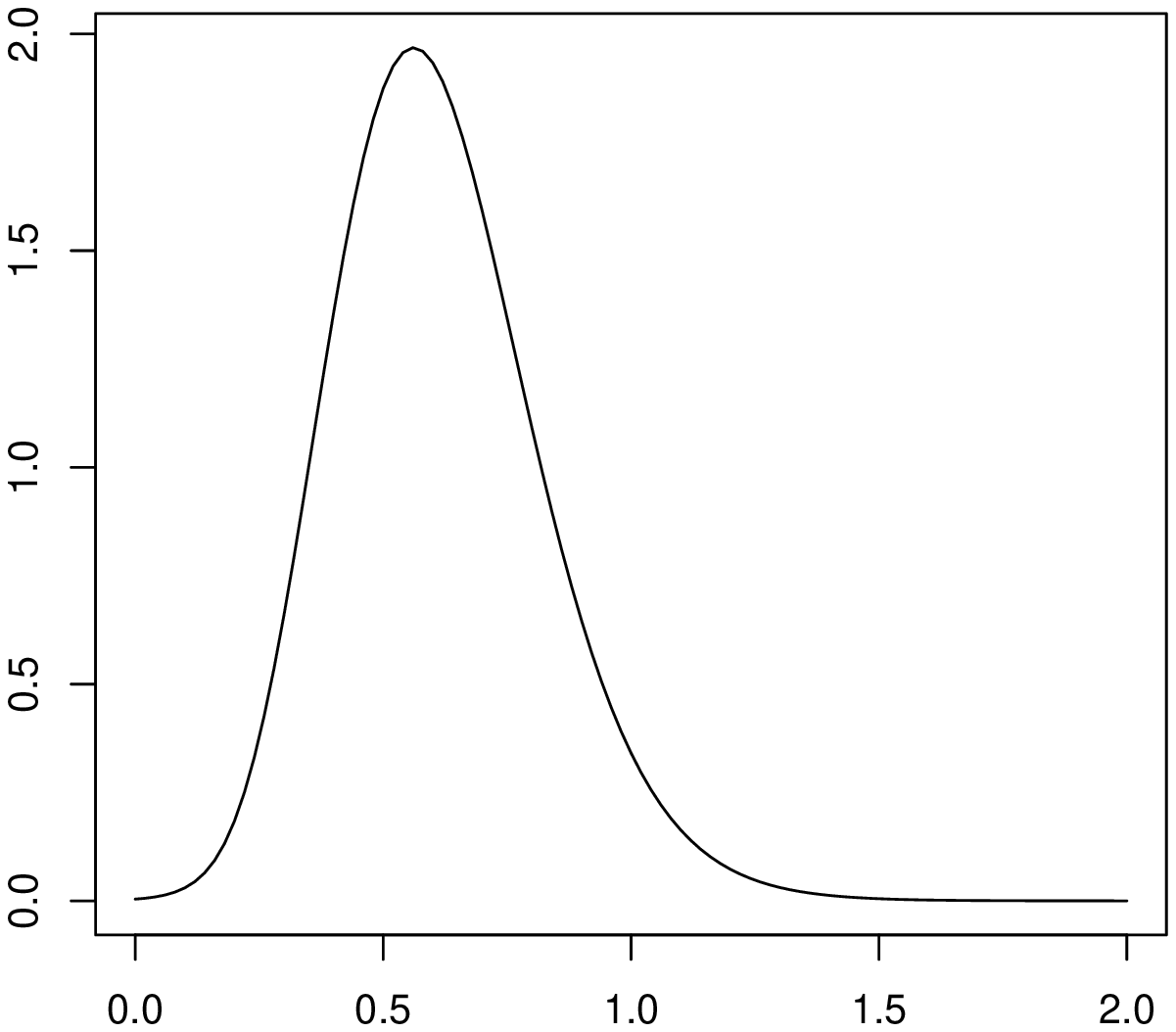}\includegraphics[scale=0.4]{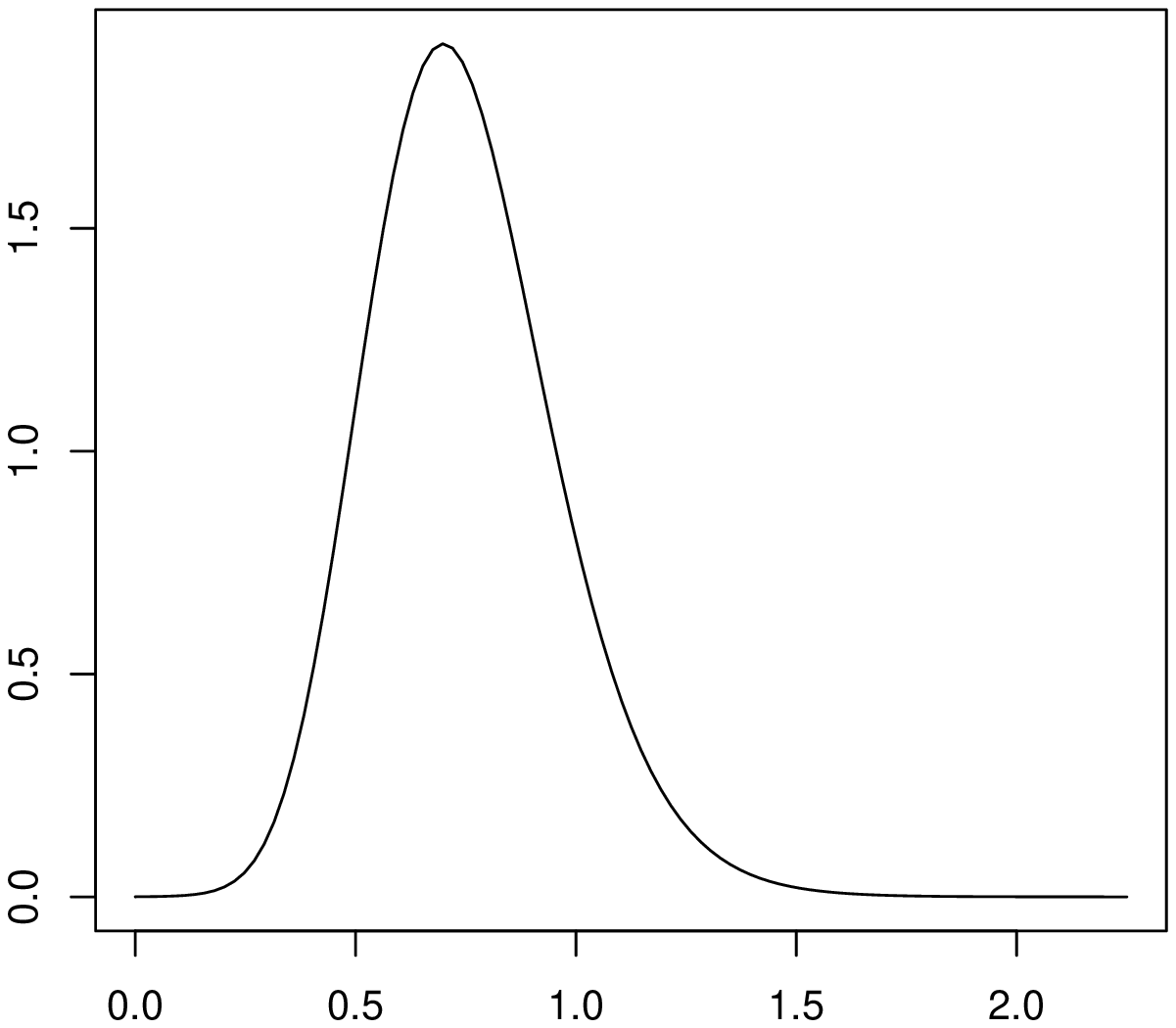}\\
%\centerline{\psfig{figure=dens1.eps,width=5cm}\psfig{figure=dens2.eps,width=5cm}\psfig{figure=dens3.eps,width=5cm}}
%\centerline{\psfig{figure=dens4.eps,width=5cm}\psfig{figure=dens5.eps,width=5cm}\psfig{figure=dens6.eps,width=5cm}}
%\centerline{\psfig{figure=dens7.eps,width=5cm}\psfig{figure=dens8.eps,width=5cm}\psfig{figure=dens9.eps,width=5cm}}
\caption{\scriptsize Predictive density of $m^{\rm th}$ repair time of a series machine with
$k$ ball bearings (in hundreds of millions revolutions) given the censored sample of Example \ref{exam1} (Scheme 1) for
 $k=1,2,3$ (rows from top to bottom) and $m=1,2,3$ (columns from left to right).}\label{f1}
\end{figure}

We shall construct equi-tailed  $95\%$ PIs, as well as $95\%$ HPD PIs, for the first
four repair times (in hundreds
of millions revolutions to failure), based on a minimal repair strategy, of a series system with
$k=1,2,3$ ball bearings, i.e. for $U_{m(k)}$,
$m=1,\ldots,4$ and $k=1,2,3$. The intervals and the corresponding widths, as well as the
point predictors $\widehat{U}_{m(k)}$, $\hat{U}^*_{m(k)}$ and
$\tilde{U}_{m(k)}$ are calculated and presented in Table
\ref{tab1}, for each of the two sampling schemes.
For example, for $k=3$ and $m=4$, the point predictor
$\widehat{U}_{4(3)}=0.9103$ means that, when a minimal repair
strategy is used, on average a series system with 3
ball bearings fails and needs to be repaired for the fourth time after
91.03 million revolutions. One can see from Table \ref{tab1} that
HPD PIs are more precise than their corresponding equi-tailed
PIs. We also provide the predictive density plots of $U_{m(k)}$, for $k=1,2,3$ and $m=1,2,3$, based on Scheme 1 in Figure \ref{f1}.

\section{Simulation study}
In this section, we wish to compare the performance of the interval and point predictors through a simulation study.
As the Rayleigh distribution belongs to location-scale family, it would be reasonable to take $\mu=0$ and $\sigma=1$ in the simulation. To examine the effect of hyper-parameters on the predictors, the informative prior is considered by setting $\xi=-1,1$ and $\tau=0.25,5$. Moreover, we consider the low-informative case by setting $\xi=0$ and $\tau=0.005$ in \eqref{eq:pdf-prior-mu-sigma}.
We take $n=20$, $r=17$ $T=2$, $k=1, 2, 3$ and $m=1, 2, 3$. The following algorithm
is used to perform the simulation:
\begin{enumerate}
\item Generate  $x_1,\cdots,x_n$ independently from Rayleigh distribution and extract the hybrid censored sample $(x_1,\cdots,x_d)$ from  $x_1,\cdots,x_n$.  

\item Generate the repair time $u_{m(k)}$ from \eqref{eq:pdfrecordexp}.

\item Obtain $95\%$ PIs as well as the point predictors $\widehat{u}_{m(k)}$, $\hat{u}^*_{m(k)}(i)$ and $\widetilde{u}_{m(k)}(i)$ based on the values of $(x_1,\cdots,x_d)$.

\item Repeat Steps 1-3 $N=10000$ times, to obtain ${u}_{m(k)}(i)$, $\widehat{u}_{m(k)}(i)$, $\hat{u}^*_{m(k)}(i)$ and $\widetilde{u}_{m(k)}(i)$, $i=1,\ldots,N$.

\item Calculate the estimated risk  (ER) of the  point predictors as follows
$$ER(\widehat{u}_{m(k)})=\frac{1}{N}\sum_{i=1}^{N}\left[\widehat{u}_{m(k)}(i)-{u}_{m(k)}(i)\right]^2,$$
$$ER(\hat{U}^*_{m(k)})=\frac{1}{N}\sum_{i=1}^{N}\left|\hat{u}^*_{m(k)}(i)-{u}_{m(k)}(i)\right|$$
and
$$ER(\widetilde{u}_{m(k)})=\frac{1}{N}\sum_{i=1}^{N}\delta\left(\widetilde{u}_{m(k)}(i)-{u}_{m(k)}(i)\right),$$
where
$$\delta(t)=\left\{\begin{array}{l l}1,& t=0\\ 0,& t\ne 0.\end{array}\right.$$

\item Calculate the average width (AW) and coverage probability (CP) of the  PIs.
\end{enumerate}

{In order to compare the Bayesian predicts with the classical (frequentist) ones, we also consider a classical method, due to Wald (1942), to obtain prediction intervals. To this end, we have to plug in the maximum likelihood estimators (MLEs) of the parameters $\mu$ and $\sigma$ based on the hybrid censored sample in the conditional (predictive) density ${f}_{U_{m(k)}|\mathbf{X}}(u|\mathbf{x})$ to estimate it and then use it to obtain the PI for $U_{m(k)}$. Using the likelihood function (\ref{eq:pdf-censorII}), the maximum likelihood estimates $\hat{\mu}$ and $\hat{\sigma}$ are obtained by solving the likelihood equations
\[\sum_{i=1}^d \frac{-1}{x_i-\mu}x_i^2+\frac{1}{\sigma}\left[\sum_{i=1}^d (x_i-\mu)+(n-d)(T_0-\mu)\right]=0\]
and
\[-\frac{d}{\sigma}+\frac{1}{2\sigma^2} \delta^*(\mu|{\mathbf{x}})=0,\]
where $\delta^*(\mu|{\mathbf{x}})$ is given in (\ref{eq:delta-star}).

Since the hybrid censored sample $\mathbf{X}$ and $U_{m(k)}$ are independent, we have $f_{U_{m(k)}|\mathbf{X}}(u|\mathbf{x})=f_{U_{m(k)}}(u)$ and thus the estimated predictive density is
obtained as
\begin{equation}
\label{eq:pdfrecordexpapp}
\widehat{f}_{U_{m(k)}}(u)=\frac{k^m (u-\hat{\mu})^{2m-1}}{\Gamma(m)\hat{\sigma}^m 2^{m-1}}\exp\left\{-\frac{k(u-\hat{\mu})^2}{2\hat{\sigma}}\right\},~~~~u>\hat{\mu}.
\end{equation}
Letting, $U_{m(k)}^*=k(U_{m(k)}-\hat{\mu})^2/\hat{\sigma}$, we have
\[\widehat{f}_{U_{m(k)}^*|\mathbf{X}}(u|\mathbf{x})=\frac{u^{m-1}e^{-u/2}}{\Gamma(m) 2^{m-1}},\]
that is the pdf of chi-square distribution with $2m$ degrees of freedom. Thus, a two-sided equi-tailed $100 (1-\alpha)\%$ Wald's prediction interval for $U_{m(k)}$, given the hybrid censored sample is
\[\left(\hat{\mu}+\sqrt{\frac{\hat{\sigma}\chi^2_{1-\alpha/2}(2m)}{k}},\hat{\mu}+\sqrt{\frac{\hat{\sigma}\chi^2_{\alpha/2}(2m)}{k}}\right),\]
where $\chi^2_{\gamma}(2m)$ denotes the upper $\gamma^{\rm th}$ quantile of
 the chi-square distribution with $2m$ degrees of freedom, i.e., $P(T>\chi^2_{\gamma}(2m))=\gamma$, with $T\sim
\chi^2(2m)$.

The simulated AWs, as well as CPs, of the Wald's PIs for $U_{m(k)}$ are given in Table 8, for $k=1,2,3$ and $m=1,2,3$. Also,
the simulated AWs and CPs of the equi-tailed Bayesian and HPD PIs for $U_{m(k)}$ are given in Tables 3-7, for different values of the hyper-parameters $\xi$ and $\tau$, $k=1,2,3$ and $m=1,2,3$.
As one can see from Tables 3-7, the ERs of the point predicts and AWs of the PIs are increasing with respect to $m$, when other parameters are kept fixed, while they are decreasing with respect to $k$. Also we observe that although HPD PIs have smaller AWs in comparison with equi-tailed PIs,  their CPs are very close and even in most cases larger than the CPs of equi-tailed PIs, revealing the superiority of HPD PIs. For  $\xi=-1$ and $\tau=5$, the largest CPs and AWs are obtained, while for $\xi=1$ and $\tau=5$, we obtain the smallest CPs and AWs. However, the ERs of all cases are almost equal. From Table 8, we see that the Wald's PIs have smaller AWs and CPs in comparison with their corresponding Bayesian PIs.

\begin{table}
\center{\caption{\textmd{The simulated AWs and CPs of the Bayesian PIs for $U_{m(k)}$, for
  $\xi=-1$ and $\tau=5$.}}}
\begin{tabular}{c c c c c c c c c c}
\hline
       &   &  Equi-tailed &   PI   &&    HPD  &   PI   &                 &                 &                 \\
 \cline{3-4}\cline{6-7}
       &$m$&    AW        &   CP   &&    AW   &   CP   &{\tiny $ER(\widehat{U}_{m(k)})$}&{\tiny $ER(\hat{U}^*_{m(k)})$}&{\tiny $ER(\widetilde{U}_{m(k)})$}\\ \hline
       &1&   3.5841  & 0.9930 && 3.4763  & 0.9875 & 0.4502 & 0.5360 &1.0000\\
$k=1$  &2&   3.8944  & 0.9890 && 3.7825  & 0.9889 & 0.5946 & 0.5970 &1.0000\\
       &3&   4.1232  & 0.9790 && 3.9736  & 0.9886 & 0.7854 & 0.6697 &1.0000\\
\hline
\hline
       &1&  2.5732  & 0.9858 &&  2.5323 & 0.9795 & 0.2396 & 0.3894 &1.0000\\
$k=2$  &2&  2.7572  & 0.9920 &&  2.6778 & 0.9876 & 0.2736 & 0.4129 &1.0000\\
       &3&  2.8789  & 0.9482 &&  2.8524 & 0.9741 & 0.5328 & 0.5600 &1.0000\\
\hline
\hline
       &1&  2.1586  & 0.9831  &&  2.1301 & 0.9750 & 0.1800 & 0.3407 &1.0000\\
$k=3$  &2&  2.2668  & 0.9909  &&  2.2410 & 0.9858 & 0.1807 & 0.3399 &1.0000\\
       &3&  2.3508  & 0.9913  &&  2.3073 & 0.9879 & 0.2029 & 0.3542 &1.0000\\
 \hline
\end{tabular}
\end{table}

\begin{table}
\center{\caption{\textmd{The simulated AWs and CPs of the Bayesian PIs for $U_{m(k)}$, for $\xi=1$ and $\tau=5$.}}}
\begin{tabular}{c c c c c c c c c c}
\hline
       &   &  Equi-tailed &   PI   &&    HPD  &   PI   &&&\\
 \cline{3-4}\cline{6-7}
       &$m$&    AW        &   CP   &&    AW   &   CP   &{\tiny $ER(\widehat{U}_{m(k)})$}&{\tiny $ER(\hat{U}^*_{m(k)})$}&{\tiny $ER(\widetilde{U}_{m(k)})$}\\ \hline
       &1&  2.4006  & 0.9036 && 2.2932 & 0.9108 &  0.4687 & 0.5423 &1.0000\\
$k=1$  &2&  2.6305  & 0.9204 && 2.5956 & 0.9278 &  0.5263 & 0.5771 &1.0000\\
       &3&  2.7528  & 0.9222 && 2.6601 & 0.9195 & 0.5698  & 0.6023 &1.0000\\
\hline
\hline
       &1&  1.7006  & 0.8899 && 1.6687 & 0.9094 & 0.2369 & 0.3865 &1.0000\\
$k=2$  &2&  1.8560  & 0.9094 && 1.8047 & 0.9172 & 0.2686 & 0.4127 &1.0000\\
       &3&  1.9450  & 0.9223 && 1.9250 & 0.9578 & 0.2848 & 0.4241 &1.0000\\
\hline
\hline
       &1&  1.4001  & 0.8804 && 1.3587 & 0.9010 & 0.1622 & 0.3195 &1.0000\\
 $k=3$ &2&  1.5176  & 0.9084 && 1.5041 & 0.9218 & 0.1772 & 0.3344 &1.0000\\
       &3&  1.5873  & 0.9187 && 1.5257 & 0.9249 & 0.1926 & 0.3471 &1.0000\\
 \hline
\end{tabular}
\end{table}
\begin{table}
\center{\caption{\textmd{The simulated AWs and CPs of the Bayesian PIs for $U_{m(k)}$, for $\xi=-1$ and $\tau=0.25$.}}}
\begin{tabular}{c c c c c c c c c c}
\hline
       &   &  Equi-tailed &   PI   &&    HPD  &   PI   &                 &                 &                 \\
 \cline{3-4}\cline{6-7}
       &$m$&    AW        &   CP   &&    AW   &   CP   &{\tiny $ER(\widehat{U}_{m(k)})$}&{\tiny $ER(\hat{U}^*_{m(k)})$}&{\tiny $ER(\widetilde{U}_{m(k)})$}\\ \hline
       &1&  2.8046  & 0.9634 &&  2.7415 & 0.9654 & 0.4619 & 0.5429 &1.0000\\
$k=1$  &2&  3.0324  & 0.9507 &&  2.9718 & 0.9561 & 0.5473 & 0.5845 &1.0000\\
       &3&  3.2103  & 0.9491 &&  3.1510 & 0.9593 & 0.6153 & 0.6157 &1.0000\\

\hline
\hline
       &1&  1.9971  & 0.9570 && 1.9859 &  0.9582 & 0.2256 & 0.3784 &1.0000\\
$k=2$  &2&  2.1387  & 0.9523 && 2.0747 &  0.9524 & 0.2657 & 0.4101 &1.0000\\
       &3&  2.2313  & 0.9519 && 2.1790 &  0.9590 & 0.2935 & 0.4252 &1.0000\\
\hline
\hline
       &1&  1.6724  & 0.9541 && 1.6606 & 0.9575 & 0.1598 & 0.3170 &1.0000\\
$k=3$  &2&  1.7622  & 0.9571 && 1.7524 & 0.9585 & 0.1764 & 0.3341 &1.0000\\
       &3&  1.8334  & 0.9489 && 1.8074 & 0.9541 & 0.1997 & 0.3532 &1.0000\\
 \hline
\end{tabular}
\end{table}
\begin{table}
\center{\caption{\textmd{The simulated AWs and CPs of the Bayesian PIs for $U_{m(k)}$, for $\xi=1$ and $\tau=0.25$.}}}
\begin{tabular}{c c c c c c c c c c}
\hline
       &   &  Equi-tailed &   PI   &&    HPD  &   PI   &                 &                 &                 \\
 \cline{3-4}\cline{6-7}
       &$m$&    AW        &   CP   &&    AW   &   CP   &{\tiny $ER(\widehat{U}_{m(k)})$}&{\tiny $ER(\hat{U}^*_{m(k)})$}&{\tiny $ER(\widetilde{U}_{m(k)})$}\\ \hline
       &1&    2.7305 & 0.9587 && 2.5025  & 0.9584  &  0.4628  & 0.5431 &1.0000\\
 $k=1$ &2&    2.9582 & 0.9428 && 2.8333  & 0.9485  &  0.5474  & 0.5795 &1.0000\\
       &3&    3.1104 & 0.9498 && 3.0390  & 0.9563  &  0.6083  & 0.6163 &1.0000\\
 \hline
\hline
       &1&   1.9378  & 0.9467 && 1.9217  & 0.9516  & 0.2313   & 0.3827 &1.0000\\
$k=2$  &2&   2.0774  & 0.9492 && 2.0129  & 0.9474  & 0.2632   & 0.4070 &1.0000\\
       &3&   2.1777  & 0.9465 && 2.1290  & 0.9525  & 0.2942   & 0.4286 &1.0000\\
\hline
\hline
       &1&  1.6144  & 0.9482  && 1.5954  & 0.9483  & 0.1582   & 0.3173 &1.0000\\
$k=3$  &2&  1.7069  & 0.9466  && 1.6974  & 0.9527  & 0.1754   & 0.3332 &1.0000\\
       &3&  1.8085  & 0.9633  && 1.7874  & 0.9678  & 0.1986   & 0.3577 &1.0000\\
\hline
\end{tabular}
\end{table}
\begin{table}
\center{\caption{\textmd{The simulated AWs and CPs of the Bayesian PIs for $U_{m(k)}$, for $\xi=0$ and $\tau=0.005$.}}}
\begin{tabular}{c c c c c c c c c c}
\hline
       &   &  Equi-tailed &   PI   &&    HPD  &   PI   &                 &                 &                 \\
 \cline{3-4}\cline{6-7}
       &$m$&    AW        &   CP   &&    AW   &   CP   &{\tiny $ER(\widehat{U}_{m(k)})$}&{\tiny $ER(\hat{U}^*_{m(k)})$}&{\tiny $ER(\widetilde{U}_{m(k)})$}\\ \hline
       &1&  2.7430  & 0.9471 &&  2.6487 & 0.9486 & 0.4661 & 0.5379 &1.0000\\
$k=1$  &2&  2.9953  & 0.9489 &&  2.9414 & 0.9539 & 0.5421 & 0.5791 &0.9999\\
       &3&  3.1791  & 0.9538 &&  3.1284 & 0.9594 & 0.6086 & 0.6118 &0.9997\\
 \hline
\hline
       &1&  1.9735  & 0.9494 && 1.9649 & 0.9548 & 0.2317 & 0.3832 &1.0000\\
$k=2$  &2&  2.1178  & 0.9452 && 2.0502 & 0.9470 & 0.2696 & 0.4107 &1.0000\\
       &3&  2.2232  & 0.9513 && 2.1658 & 0.9559 & 0.2960 & 0.4293 &1.0000\\
\hline
\hline
       &1& 1.6489   & 0.9493  && 1.6334 & 0.9525 & 0.1590 & 0.3173 &1.0000\\
$k=3$  &2& 1.7390   & 0.9526  && 1.7301 & 0.9571 & 0.1733 & 0.3313 &1.0000\\
       &3& 1.8103   & 0.9503  && 1.7897 & 0.9538 & 0.1974 & 0.3482 &1.0000\\
 \hline
\end{tabular}
\end{table}

 {

\begin{table}
\center{\caption{\textmd{The simulated AWs and CPs of the Wald's PIs for $U_{m(k)}$.}}}
\begin{tabular}{c c c c }
\hline
       &$m$&    AW        &   CP   \\
        \hline
       &1&  2.3382  & 0.9059 \\
$k=1$  &2&  2.4778  & 0.9091 \\
       &3&  2.5293  & 0.8943 \\
 \hline
\hline
       &1&  1.6496  & 0.8964 \\
$k=2$  &2&  1.7597  & 0.9027 \\
       &3&  1.7807  & 0.8978 \\
\hline
\hline
       &1& 1.3463   & 0.8801  \\
$k=3$  &2& 1.4294   & 0.8994  \\
       &3& 1.4734   & 0.9022  \\
 \hline
\end{tabular}
\end{table}

}
\section{Concluding Remarks}
In this paper, we have studied the marginal Bayesian
prediction of repair times of a series system based on a minimal
repair strategy, using the information contained in an observed
hybrid censored sample of the lifetimes of the components of the
system. The results can be extended to non-parametric prediction 
procedures and also the joint prediction of
repair times. The problem of predicting repair times
of other types of coherent systems, such as parallel systems and $k$-out-of-$n$ systems, on the basis of a censored
sample of their components is of interest. Work on these
problems is under progress and we hope to report the new results
in the near future.

\section*{Acknowledgements}
We would like to thank the associate editor and the anonymous referees for
their valuable comments and suggestions on an earlier version of this manuscript, which resulted in this improved version. { We also thank Prof. N. Balakrishnan for his useful comments.}

 \section*{References}
\begin{enumerate}
\renewcommand{\baselinestretch}{1}

\item \vskip  -2mm

Abdel-Hamid, A.H. and AL-Hussaini, E.K. (2014). Bayesian prediction for Type-II progressive-censored data from the Rayleigh distribution under progressive-stress model. \textit{Journal of Statistical Computation and Simulation}., \textbf{84}, 1297--1312.

\item \vskip  -2mm

Ahmadi, J. and Balakrishnan, N. (2010). Prediction of order
statistics and record values from two independent sequences,
\textit{Statistics}, \textbf{44}, 417--430.

\item \vskip  -2mm

Ahmadi, J. and MirMostafaee, S. M. T. K. (2009).
Prediction intervals for future records and order statistics coming  from  two
parameter exponential distribution, \textit{Statistics and Probability Letters}, \textbf{79}, 977--983.

\item \vskip  -2mm

Ahmadi, J. and  MirMostafaee, S. M. T. K. and Balakrishnan, N.
(2010). Nonparametric prediction intervals for future record
intervals based on order statistics, \textit{Statistics and Probability Letters}, \textbf{80}, 1663--1672.

\item \vskip  -2mm

Ali Mousa, M. A. M., Al-Sagheer, S. A. (2006). Statistical inference for the Rayleigh model
based on progressively Type-II censored data. \textit{Statistics}, \textbf{40}, 149--157.

\item \vskip  -2mm 

Arnold, B. C.,  Balakrishnan, N., and  Nagaraja, H.
N. (1998). {\it Records}, John Wiley \& Sons, New York.

\item \vskip  -2mm

Balakrishnan, N. (1989). Approximate MLE of the scale parameter of the Rayleigh distribution with censoring. \textit{IEEE Transactions on Reliability}, \textbf{38}, 355--357.

\item \vskip  -2mm

Balakrishnan, N. and Kundu, D. (2012). Hybrid
censoring: models, inferential results  and applications (with
discussion). {\it Computational Statistics and Data Analysis},
\textbf{57}, 166--209.

\item \vskip  -2mm

Barlow, R. E. and Hunter, L. (1960). Optimum preventive maintenance policies. {\it Operations Research},
\textbf{8}, 90--100.

\item \vskip  -2mm

Childs, A., Chandrasekhar, B., Balakrishnan, N. and Kundu,
D. (2003). Exact likelihood inference based on type-I and type-II
hybrid censored samples from the exponential distribution.
\textit{Annals of the Institute of Statistical Mathematics},
\textbf{55}, 319--30.

\item \vskip  -2mm

Dey, S. and Das, M.K. (2007). A Note on prediction interval for a Rayleigh distribution: Bayesian approach. {\it American Journal of Mathematical and Management Sciences}, \textbf{27}, 43--48.

\item \vskip  -2mm

Dey, S. and Dey, T. (2012). Bayesian estimation and prediction intervals for a Rayleigh distribution under a conjugate prior. \textit{Journal of Statistical Computation and Simulation}., \textbf{82}, 1651--1660.

\item \vskip  -2mm

Dey, S. and Dey, T. (2014). Statistical inference for the Rayleigh distribution under Type II progressive censoring with binomial removal. {\it Applied Mathematical Modelling}, {\bf 38}, 974--982.

\item \vskip  -2mm

Draper, N. and Guttman, I. (1987). Bayesian analysis
of hybrid life tests with exponential failure times, \textit{Annals
of the Institute of Statistical Mathematics},  \textbf{39},
219--225.
\item \vskip  -2mm

Dyer, D.D. and Whisenand, C.W. (1973). Best linear unbiased estimator of the parameter of the Rayleigh distribution-Part I: small sample theory for censored order statistics. \textit{IEEE Transactions on Reliability}, \textbf{22}, 27--34.

\item \vskip -2mm

Ebrahimi, N. (1992). Prediction intervals for future
failures in exponential distribution under hybrid censoring,
\textit{IEEE Transactions on Reliability}, \textbf{41}, 127--132.

\item \vskip  -2mm

Epstein, B. (1954). Truncated life tests in the exponential case.
\textit{Annals of Mathematical Statistics},
\textbf{25}, 555--564.

\item \vskip  -2mm

Fernandez, A. J. (2000). Bayesian inference from Type-II doubly censored Rayleigh data.
\textit{Statistics and Probability Letters}, \textbf{48}, 393--399.

\item \vskip  -2mm

Gelman, A., Carlin, J. B., Stern, H. S. and Rubin, D. B. (2004). {\it Bayesian Data Analysis, Second Edition}, Chapman \& Hall/CRC, New York, USA. 

\item \vskip  -2mm

Gross, A.J. and Clark, V.A. (1976). {\it Survival Distributions, Reliability Applications in Biomedical Sciences (Probability and Mathematical Statistics)}, John Wiley \& Sons, New York.

\item \vskip  -2mm

Howlader, H.A. (1985). HPD prediction intervals for the Rayleigh distribution. \textit{IEEE Transactions on Reliability}, \textbf{34}, 121--123.

\item \vskip  -2mm

Howlader, H.A. and Hossain, A. (1995). On bayesian estimation and prediction from rayleigh based on Type II censored data. \textit{Communications in Statistics- Theory and Methods}, \textbf{24}, 2251--2259.

\item \vskip  -2mm

Johnson, N.L., Kotz, S. and Balakrishnan, N. (1994). {\it Continuous Univariate Distributions}, Vol. 1, 2nd Ed., John Wiley \& Sons, New York.

\item \vskip  -2mm

Khan, H. M. R., Provost, S. B. and Singh, A. (2010). Predictive inference from a two-Parameter Rayleigh
life model given a doubly censored sample. \textit{Communications in Statistics- Theory and Methods}, \textbf{39}, 1237--1246.

\item \vskip  -2mm

Kundu, D. (2007). On hybrid censoring Weibull
distribution., \textit{Journal of Statistical Planning and
Inference}, \textbf{137}, 2127--2142.

\item \vskip  -2mm

Kundu, D. and Banerjee, A. (2008). Inference Based
on Type-II Hybrid Censored Data From a Weibull Distribution.
\textit{IEEE Transactions on Reliability}, \textbf{57(2)}, 369--378.

\item \vskip  -2mm

Kundu, D. and Gupta, R. D.  (1988). Hybrid
censoring schemes with exponential failure distribution.
\textit{Communications in Statistics- Theory and Methods},
\textbf{27}, 3065--3083.

\item \vskip  -2mm

Lawless, J.F. (2003). {\it Statistical Models and Methods for Lifetime Data}, John Wiley \& Sons, New York.

\item \vskip  -2mm

Leiblein, J. and Zelen, M. (1956). Statistical investigation of the fatigue life of deep-groove ball
bearings, \textit{Journal of Research of the National Bureau of Standards}, \textbf{57}, 273--316.

\item \vskip  -2mm

MirMostafaee, S. M. T. K. and Ahmadi, J. (2011). Point prediction
of future order statistics  from exponential distribution,
\textit{Statistics and Probability Letters}, \textbf{81}, 360--370.

\item \vskip -2mm

Nevzorov,  V. (2001). {\it Records: Mathematical Theory}, Translation of Mathematical
Monographs No. 194, American
Mathematical Society, Providence, Rhode Island, USA.

\item \vskip -2mm

Polovko, A. M. (1968). {\it Fundamentals of Reliability Theory}, Academic Press, San Diego.

\item \vskip  -2mm

Raqab, M. Z. and Madi, M. T. (2002). Bayesian prediction of the total time on test using doubly censored Rayleigh
data. \textit{Journal of Statistical Computation and
Simulation}., \textbf{72}, 781--789.

\item \vskip  -2mm

Rayleigh, J.W.S. (1880). On the resultant of a large number of vibrations of the same pitch and of arbitrary phase. {\it The London, Edinburgh, and Dublin Philosophical Magazine and Journal of Science}, 5th Series, \textbf{10}, 73--78.

\item \vskip  -2mm

Rayleigh, J.W.S. (1919). On the problem of random variations, and of random flights in one, two or three dimensions. {\it The London, Edinburgh, and Dublin Philosophical Magazine and Journal of Science}, 6th Series, \textbf{37}, 321--347.

\item \vskip  -2mm

Soliman, A. A. and Al-Aboud, F. M. (2008). Bayesian inference using record values
from Rayleigh model with application. \textit{European Journal of Operational Research}, \textbf{185}, 659--672.

\item \vskip  -2mm

{
Wald, A.  (1942). Setting of tolerance limits when the sample is large. \textit{Annals of Mathematical Statistics}, \textbf{13}, 389--399.

}
\end{enumerate}

\end{document}